\documentclass{jnmp}
\usepackage{amsmath}

\begin{document}
\setcounter{page}{229}
\renewcommand{\evenhead}{A~Sergeev}
\renewcommand{\oddhead}{Enveloping Superalgebra $U({\mathfrak{osp}}(1|2))$ and Orthogonal
Polynomials}

\thispagestyle{empty}


\FirstPageHead{8}{2}{2001}
{\pageref{sergeev-firstpage}--\pageref{sergeev-lastpage}}{Article}

\copyrightnote{2001}{A~Sergeev}

\Name{Enveloping Superalgebra ${U({\mathfrak{osp}}(1|2))}$
and\\ Orthogonal Polynomials
in Discrete Indeterminate}\label{sergeev-firstpage}

\Author{A~SERGEEV}

\Address{(Correspondence) D Leites, Department of Mathematics, University of Stockholm\\
Roslagsv.  101, Kr\"aftriket hus 6, SE-106 91, Stockholm, Sweden\\
E-mail: mleites@matematik.su.se\\[2mm]
On leave of absence from Balakovo Institute of Technique of Technology and Control,\\
Branch of Saratov Technical University, Balakovo, Saratov Region, Russia}

\Date{Received May 6, 2000; Accepted  November 18, 2000}

\begin{abstract}
\noindent
Let $A$ be an associative simple (central)
superalgebra over ${\mathbb C}$ and $L$ an invariant linear functional on it
(trace).  Let $a\mapsto a^t$ be an antiautomorphism of $A$ such that
$(a^t)^ t=(-1)^{p(a)}a$, where $p(a)$ is the parity of $a$, and let
$L(a^t)=L(a)$.  Then $A$ admits a nondegenerate supersymmetric
invariant bilinear form $\langle a, b\rangle=L(ab^t)$.  For
$A=U({\mathfrak{sl}}(2))/{\mathfrak{m}}$, where ${\mathfrak{m}}$ is any maximal ideal of $U({\mathfrak{sl}}(2))$,
Leites and I have constructed orthogonal basis in $A$ whose elements
turned out to be, essentially, Chebyshev (Hahn) polynomials in one
discrete variable.
Here I take $A=U({\mathfrak{osp}}(1|2))/{\mathfrak{m}}$ for any maximal ideal ${\mathfrak{m}}$ and apply
a similar procedure.  As a result we obtain either Hahn polynomials
over ${\mathbb C}[\tau]$, where $\tau^2\in{\mathbb C}$, or a particular case of
Meixner polynomials, or --- when $A=\mbox{Mat}(n+1|n)$ --- dual Hahn
polynomials of even degree, or their (hopefully, new) analogs of odd
degree.  Observe that the nondegenerate bilinear forms we consider
for orthogonality are, as a rule, not sign definite.
\end{abstract}

\renewcommand{\theequation}{\thesection.\arabic{equation}}
\setcounter{equation}{0}

\section{Introduction}

Classically, orthogonal polynomials were considered with respect to a
sign definite bilinear form.  Lately we encounter the growth of
interest to the study of orthogonal polynomials relative an arbitrary
(but still symmetric and nondegenerate) form, cf.~\cite{sergeev:KLP,sergeev:M,sergeev:MBP}
and refe\-ren\-ces therein.  In these approaches,
however, the bilinear forms are introduced ``by hands" and the
differential or difference equations the orthogonal polynomials
satisfy are of high degree.  We would like to point out that traces
and supertraces on associative algebras and superalgebras are natural
sources of bilinear symmetric forms which are seldom sign-definite.
The Lie structure on the algebras obtained from these associative
algebras and superalgebras is more adapted to the study of orthogonal
polynomials.  In particular, the eigenvalue problem for the Casimir
operator --- the quadratic element of the center with respect to the
Lie structure --- naturally provides with a 2nd degree difference
equation for the polynomials orthogonal relative the above
(super)traces.

Let ${\mathfrak{sl}}(2)$ be represented as $\mbox{Span}(Y, H, X)$
subject to relations
\begin{equation}
[X, Y]=H, \qquad [H, Y]=-2Y, \qquad [H, X]=2X.\label{1.1.1}
\end{equation}
The quadratic Casimir operator of ${\mathfrak{sl}}(2)$
\begin{equation}
\Omega =2YX+\frac12 H^2+ H\label{1.1.2}
\end{equation}
lies in the center of $U({\mathfrak{sl}}(2))$.  Let $I_{\lambda}$ be the
two-sided ideal in the associative algebra $U({\mathfrak{sl}}(2))$ generated by
$\Omega-\frac12(\lambda^{2}-1)$.  It turns out that the associative
algebra $\tilde{\mathfrak{A}}_{\lambda}=U({\mathfrak{sl}}(2))/I_{\lambda}$ is simple for
$\lambda\not\in {\mathbb{Z}}\setminus\{0\}$, otherwise $\tilde{\mathfrak{A}}_{\lambda}$
contains an ideal such that the quotient is isomorphic to the matrix
algebra $\mbox{Mat}(|\lambda|)$.  Set (\cite{sergeev:Di})
\begin{equation}
{\mathfrak{A}}_{\lambda}=\left\{\begin{array}{ll}
\tilde{\mathfrak{A}}_{\lambda}, &\mbox{if} \ \
\lambda\not\in {\mathbb{Z}}\setminus\{0\}, \vspace{1mm}\\
\mbox{Mat}(|\lambda|), \quad & \mbox{otherwise}.
\end{array}\right.\label{1.1.3}
\end{equation}
Clearly, ${\mathfrak{A}}_{-\lambda}\simeq{\mathfrak{A}}_\lambda$.  As associative algebra,
${\mathfrak{A}}_{\lambda}$ is generated by $X$, $Y$, and $H$ subject to relations
\begin{equation}
XY=\frac14\left(\lambda^2-(H-1)^2\right) \label{1.1.4}
\end{equation}
and one more relation for integer values of $\lambda$:
\begin{equation}
X^{|\lambda|}=0 \quad \mbox{if}\ \ \lambda\in
{\mathbb{Z}}\setminus\{0\}.\label{1.1.5}
\end{equation}
It is also known that ${\mathfrak{A}}_\lambda$ possesses an antiautomorphism
$u\mapsto u^t$ given on generators by the formula
\begin{equation}
X^t= Y, \qquad Y^t=X, \qquad H^t=H.\label{1.1.6}
\end{equation}
In \cite{sergeev:LS} we have shown that on ${\mathfrak{A}}_\lambda$ there exists a
unique, up to a constant factor, nontrivial linear functional~$L$,
which for positive integer $\lambda$'s is the usual trace and which
satisfies $L\left(u^t\right)= L(u)$.  By means of this functional we define an
invariant symmetric bilinear form $\langle\cdot, \cdot\rangle$ on
${\mathfrak{A}}_\lambda$, by setting $\langle u, v\rangle=L\left(uv^t\right)$.  The form
$\langle\cdot, \cdot\rangle$ is nondegerate and symmetric.

Now, consider ${\mathfrak{A}}_\lambda$ as an ${\mathfrak{sl}}(2)$-module with respect to the
adjoint representation.  We have
\begin{equation}
{\mathfrak{A}}_{\lambda}=\left\{\begin{array}{ll}\sum\limits_{i=0}^\infty
L^{2i}, &\mbox{if}\ \  \lambda\not\in {\mathbb{Z}}\setminus\{0\},
\vspace{2mm}\\
\sum\limits_{i=0}^{|\lambda|-1} L^{2i}, \quad &\mbox{otherwise}.
\end{array}\right.\label{1.1.7}
\end{equation}
where $L^{2i}$ is the irreducible finite dimensional ${\mathfrak{sl}}(2)$-module
with highest weight $2i$ (cf.~\cite{sergeev:F}).  Clearly, $H$ arranges a
${\mathbb{Z}}$-grading on ${\mathfrak{A}}_{\lambda}$, namely
\begin{equation}
({\mathfrak{A}}_{\lambda})_{i}=\{u\in {\mathfrak{A}}_{\lambda}\mid [H, u]=2iu \ \mbox{for}
\ i\in{\mathbb{Z}}\}.\label{1.1.8}
\end{equation}

For any $f, g\in{\mathbb C}[H]$ and $i\geq 0$  set
\begin{equation}
\langle f, g\rangle_{i}= \langle f(H)X^i, g(H)X^i\rangle\qquad
\mbox{and}\qquad \langle f, g\rangle_{-i}= \langle f(H)Y^i, g(H)Y^i\rangle.
\label{1.1.9}
\end{equation}
Denote:
\begin{equation}
T_{i}(H)=\frac14\left(n^{2}-(H+2i-1)^2\right)\qquad \mbox{and}\qquad
\alpha_{i}=n-2i+1.\label{1.1.10}
\end{equation}
Set
\begin{equation}
f_{ki}(H)X^i=(\mbox{ad}\; Y)^k\left(X^{k+i}\right) ;  \label{1.1.11}
\end{equation}
further set
\begin{equation}
(\Delta_{2}f)(H)=f(H+2)-f(H) \qquad \text{and}\qquad
(\nabla_{2}f)(H)=f(H)-f(H-2). \label{1.1.12}
\end{equation}

\noindent
{\bf Theorem 1.1 \cite{sergeev:LS}.} {\it

\vspace{-2mm}

\begin{enumerate}
\topsep0mm
\partopsep0mm
\parsep0mm
\itemsep0mm

\item[1)] $\langle ({\mathfrak{A}}_{\lambda})_{i},
({\mathfrak{A}}_{\lambda})_{j}\rangle=0$ for $i\neq j$.

\item[2)] For $i\geq 0$ the polynomials $f_{ki}(H)$ are of degree $k$,
they are orthogonal relative to the form $\langle \cdot ,
\cdot\rangle_{i}$.

\item[3)] For $i\geq 0$ the polynomials $f_{ki}(-H)$ are of degree $k$,
they are orthogonal relative to the form $\langle \cdot ,
\cdot\rangle_{-i}$.

\item[4)] The polynomials $f_{ki}(H)$ satisfy the following difference
equation:
\begin{equation}
\displaystyle \frac14\left((H+1)^2\!-\!\lambda^2\right)\!
 (\Delta_{2}f)\!-\!\frac14\left((H\!-\!2i\!-\!1)^2\!-\!\lambda^2\right)
(\nabla_{2}f)
\!=\!k(k+2i+1)f.
\label{1.2.1}
\end{equation}

\item[5)] Explicitly we have
\begin{equation}\arraycolsep=0em
\begin{array}{l}
f_{kl}(H)={}_{3}F_{2}\left(\left.\begin{array}{c} l-k, \; l+k+1, \;
\frac12(1-n-H)\vspace{1mm}\\ l+1, \;
l+1-n\end{array}\, \right| \, 1\right)
\vspace{2mm}\\
\phantom{f_{kl}(H)={}} {}\times T_{0}((\alpha)_{l+1})\cdots
T_{0}((\alpha)_{l+k}),
\end{array}\label{1.2.2}
\end{equation}
where
\[
{}_{3}F_{2}\left(\left.\begin{array}{c}\alpha_{1}, \; \alpha_{2}, \; \alpha_{3}\vspace{1mm}\\
\beta_{1}, \; \beta_{2}\end{array}\right|\;
z\right)=\mathop{\sum}\limits_{i=0}^\infty
\displaystyle\frac{(\alpha_{1})_{i}(\alpha_{2})_{i}
(\alpha_{3})_{i}}{(\beta_{1})_{i}(\beta_{2})_{i}}\,
\displaystyle\frac{z^i}{i!}
\]
is a generalized hypergeometric function, $(\alpha)_{0}=1$  and
$(\alpha)_{i}=\alpha(\alpha+1)\cdots(\alpha+i-1)$ for $i>0$.
\end{enumerate}}

{\bf Our goal is to generalize this theorem by replacing
$\pbf{{\mathfrak{sl}}(2)}$ with $\pbf{{\mathfrak{osp}}(1|2)}$.}
 The main result obtained is the union of Theorems 2.2, 2.3 and 2.4.

\setcounter{equation}{0}

\section{Preliminaries and main result}

We select the following basis in ${\mathfrak{osp}}(1|2)\subset
{\mathfrak{sl}}(\bar 0|\bar 1|\bar 0)$:
\[
Y=\left(\begin{array}{rcc}
0&0&0\\
0&0&0\\
-1&0&0\end{array}\right), \qquad
H=\left(\begin{array}{ccr}
1&0&0\\0&0&0\\
0&  0&-1 \end{array}\right), \qquad X= \left(\begin{array}{ccc} 0& 0&1\\
0&0&0\\
0&0&0
\end{array}\right)
\]
\[
F=\left(\begin{array}{ccc} 0&0&0\\
1&0&0\\
0&1&0
\end{array}\right), \qquad G=\left(\begin{array}{ccr}
0&1&0\\
0&0&-1\\
0&0&0 \end{array}\right).
\]
The defining relations (we give only the ones with nonzero values in
the right hand side) are
\begin{equation}\arraycolsep=0em
\begin{array}{l}
{}[H, F]=-F, \qquad [H, G]=G, \qquad [G, F]=H,
\vspace{2mm}\\
{}[G, G]=2X, \qquad [F, F]=-2Y.
\end{array}\label{2.1.1'}
\end{equation}
For convenience we add also the following corollaries
\begin{equation}\arraycolsep=0em
\begin{array}{l}
{}[H, X]=2X, \qquad [H, Y]=-2Y,
\vspace{2mm}\\
{}[X, Y]=H, \qquad [G, Y]=F, \qquad [F, X]=G.
\end{array}\label{2.1.1''}
\end{equation}
These relations immediately imply that $U({\mathfrak{osp}}(1|2))$ is generated, as
associative superalgebra, by $F$, $H$, and $G$.  Set
\begin{equation}\arraycolsep=0em
\begin{array}{l}
\displaystyle \tau=2GF-H+\frac12=H+\frac12-2FG, \vspace{2mm}\\
\Omega=H^2-H+4XY-2GF=H^2+H+4XY-2FG, \vspace{2mm}\\
\omega =H^2-2H+4XY= H^2+2H+4YX.
\end{array}\label{2.1.2}
\end{equation}

\noindent
{\bf Lemma 2.1.} {\it
\vspace{-2mm}

\begin{enumerate}
\topsep0mm
\partopsep0mm
\parsep0mm
\itemsep0mm
\item[i)]  The element $\tau$ belongs to the
supercenter of $U({\mathfrak{osp}}(1|2))$, i.e., it commutes with the even
elements and anticommutes with the odd ones.

\item[ii)] The centralizer of the Cartan subalgebra of
${\mathfrak{osp}}(1|2)$ in $U({\mathfrak{osp}}(1|2))$ is generated
by~$H$ and $\tau$.

\item[iii)] $\Omega= \tau^2-\frac14$ is the quadratic Casimir element
of $U({\mathfrak{osp}}(1|2))$.

\item[iv)] $\omega = \tau^2+\tau-\frac34$ is the quadratic Casimir
element of $U({\mathfrak{sl}}(2))$.
\end{enumerate}}

\noindent {\bf Proof.} i) (the proof of the fact that $\tau G+G\tau=0$
is similar): \[\arraycolsep=0em
\begin{array}{l}
   \displaystyle \tau
   F+F\tau=\left(2GF-H+\frac12\right)F+F\left(-2FG+H+\frac12\right)
   \vspace{2mm}\\
\phantom{\tau F+F\tau} {}= FH-HF+F+2GY-2YG=2F-2[G, Y]=0.
\end{array}
\]

ii) It is easy to verify that
\begin{equation}\arraycolsep=0em
\begin{array}{l}
\displaystyle G^nF^n=\left(\frac12\right)^n\left(\tau+H-\frac12\right)\left(-\tau+H-
\frac32\right)
\vspace{3mm}\\
\displaystyle \phantom{G^nF^n=} {}\times \left(\tau+H-\frac52\right)\cdots
\left ((-1)^n\tau+H-\frac{2n-1}{2}\right )\label{2.2.1}
\end{array}
\end{equation}
Now observe that any element of the centralizer is a linear
combination of the elements $G^nF^n$ for $n\geq 0$.

Headings iii) and iv) are subject to a similar direct verification. \hfill \rule{3mm}{3mm}

\medskip

{\bf  A theorem of Pinczon.} Pinczon~\cite{sergeev:Pi} described the
maximal two-sided ideals of $U({\mathfrak{osp}}(1|2))$.  Let us formulate his
results in a form convenient to us.

\newpage

\noindent
{\bf Theorem 2.1.} {\it

\vspace{-2mm}

\begin{enumerate}
\topsep0mm
\partopsep0mm
\parsep0mm
\itemsep0mm
\item[A)] Every maximal two-sided ideal of $U({\mathfrak{osp}}(1|2))$ is of the form
${\mathfrak{m}}_\lambda$, where ${\mathfrak{m}}_\lambda$ is:\\
    generated by $\tau^2-\left(\lambda+\frac12\right) ^2$ for
$\lambda\not\in{\mathbb{Z}}_{\geq 0}\cup \left\{-\frac12\right\}$;\\
 the kernel of the finite dimensional representation with highest
weight $\lambda$ for $\lambda\in{\mathbb{Z}}_{\geq 0}$;\\
generated by $\tau$ for $\lambda=-\frac12$.
\item[B)] Let ${\mathfrak{B}}_\lambda=U({\mathfrak{osp}}(1|2))/{\mathfrak{m}}_\lambda$.  Then
 \begin{enumerate}
\item[i)] If $\lambda\not\in{\mathbb{Z}}_{\geq 0}\cup\{-\frac12\}$, then
${\mathfrak{B}}_\lambda$ is generated by $G$, $H$, and $F$ subject to the
relations
\begin{equation}\arraycolsep=0em
\begin{array}{l}
\displaystyle
{}[H, G]=G, \qquad [H, F]=-F;\qquad [G, F]=H;
\vspace{2mm}\\
\displaystyle  \tau=H+\frac12-2FG, \qquad
\tau^2=\left(\lambda+\frac12\right)^2.
\end{array}\label{2.3.1}
\end{equation}
\item[ii)] The superalgebra ${\mathfrak{B}}_{-\frac12}$ is isomorphic to the Weil
algebra $A_1=\mbox{\rm Diff}\, (1)={\mathbb C}[P, Q]$ considered as superalgebra when
generators are considered to be odd (recall that the defining
relations in $\mbox{\rm Diff}\, (1)$ are $PQ-QP=1$.
 \item[iii)] ${\mathfrak{B}}_\lambda\simeq \mbox{\rm Mat}(\lambda+1|\lambda)$ for
$\lambda\in{\mathbb{Z}}_{\geq 0}$.
\end{enumerate}
\end{enumerate}}

{\bf On the structure of $\pbf{{\mathfrak{B}}_\lambda}$.} Recall that an {\it
antiautomorphism of super}algebra $A$ is an even linear map $a\mapsto
a^t$ for $a\in A$ such that $(ab)^t=(-1)^{p(a)p(b)} b^ta^t$.  Define
an antiautomorphism of $U({\mathfrak{osp}}(1|2))$ by setting $H^t=H$, $F^t=G$,
$G^t=-F$.  Clearly, this antiautomorphism induces an antiautomorphism
of ${\mathfrak{B}}_\lambda$ for every $\lambda$.

Later on, I will show that on ${\mathfrak{B}}_\lambda$ exists a unique, up to a
scalar factor, nontrivial invariant linear functional --- the
supertrace $\mbox{str}$.  So the form $\langle u, v\rangle=\mbox{str}\left(uv^t\right)$
determines an invariant supersymmetric bilinear form on ${\mathfrak{B}}_\lambda$.
The lack of nonzero two-sided ideals guarantees the non-degeneracy of the
form.

For $\lambda\not\in{\mathbb{Z}}_{+}\cup\left\{-\frac12\right\}$ the algebra
${\mathfrak{B}}_\lambda$ possesses a ${\mathbb{Z}}$-grading
of the form ${\mathfrak{B}}_\lambda=
\mathop{\oplus}\limits_{i\in{\mathbb{Z}}}({\mathfrak{B}}_\lambda)_{i}$, where
\[
({\mathfrak{B}}_\lambda)_{i}=\{u\in {\mathfrak{B}}_\lambda\mid [H, u]=iu\ \ \text{for}
\ \  i\in{\mathbb{Z}}\}.
\]
For $f, g\in{\mathbb C}[H, \tau]$ and $i\geq 0$ set
\begin{equation}
\langle f, g\rangle_{i}= \langle fG^i, gG^i\rangle\qquad\text{and}\qquad
\langle f, g\rangle_{-i}= \langle fF^i, gF^i\rangle.\label{2.4.1}
\end{equation}
Recall that $\tau^2=\left(\lambda+\frac12\right)^2$ and introduce a
${\mathbb{Z}}$-grading in ${\mathbb C}[H, \tau]$ by setting $\deg H=2$, $\deg
\tau=1$.  Now it is not difficult to verify that for any $i\in{\mathbb{Z}}$
there exists a basis $\{f_k\}_{k\geq 0}$ such that $\deg f_k=k$ and
$\langle f_k, f\rangle=0$ if $\deg f<k$.  In what follows the basis
elements with such properties will be called {\it orthogonal
polynomials} in $H$ and $\tau$.

Set also
\begin{equation}
\Delta(H)=\nabla(H)=1, \qquad \Delta(\tau)=-2\tau, \qquad
\nabla(\tau)=2\tau\label{2.4.2}
\end{equation}
and extend the action of the operators $\Delta$ and $\nabla$ onto
${\mathbb C}[H, \tau]$ by setting
\begin{equation}\arraycolsep=0em
\begin{array}{l}
\Delta(fg)=\Delta(f)\cdot g+ f\cdot \Delta(g)+\Delta(f)\Delta(g),
\vspace{2mm}\\
\nabla(fg)=\nabla(f)\cdot g+ f\cdot \nabla(g)+\nabla(f)\nabla(g).
\end{array}\label{2.4.3}
\end{equation}
Besides, set
\[
\Delta_{2}(\tau)=\nabla_{2}(\tau)=0
\]
and extent the action of the operators $\Delta_{2}$ and $\nabla_{2}$
onto ${\mathbb C}[H, \tau]$ by formulas similar to~(\ref{2.4.3}).

Define now polynomials $f_{k, j}$ for $j\geq 0$ by setting
\begin{equation}
f_{k, 2i}G^{2i}=\left\{\begin{array}{ll} (\mbox{ad}\;F)^k\left(G^{2i+k}\right)&\text{for $k$ even},
\vspace{1mm}\\
(\mbox{ad}\; F)^k\left(G^{2i+k}\tau\right)\quad &\text{for $k$ odd}
\end{array}\right.\label{2.4.4}
\end{equation}
and
\begin{equation}
f_{k, 2i+1}G^{2i+1}=\left\{\begin{array}{ll}
(\mbox{ad}\; F)^{k+1}\left(G^{2i+k+2}\right)\quad &\text{for $k$ even}, \vspace{1mm}\\
(\mbox{ad}\; F)^{k+1}\left(G^{2i+k}\tau\right)&\text{for $k$ odd}.
\end{array}\right.\label{2.4.5}
\end{equation}

\noindent
{\bf Theorem 2.2.} {\it
\vspace{-2mm}

\begin{enumerate}
\topsep0mm
\partopsep0mm
\parsep0mm
\itemsep0mm
\item[1)] $\langle ({\mathfrak{B}}_\lambda)_{i},
({\mathfrak{B}}_\lambda)_{j}\rangle=0$ for $i\neq j$.

\item[2)] Polynomials $f_{k, j}(H)$ are orthogonal relative to the form
$\langle \cdot, \cdot \rangle_{i}$; the polynomial $f_{k, j}(H)$ is a
degree $k$ polynomial in $H$ and $\tau$.

\item[3)] Polynomial $f_{k, 2i}(H)$ satisfies the difference equation
\[
\left(H-\tau+\frac12\right)\Delta\nabla f+2(H-i)\nabla
    f+\left(2i+\frac12\right)f=(-1)^{p(k)}\left(2i+k+\frac12\right)f.
\]

\item[4)] Polynomial $f_{k, 2i+1}(H)$ satisfies the difference equation
\[\arraycolsep=0em
\begin{array}{l}
\displaystyle \left(H-\tau+\frac12\right)\Delta\nabla f+2\left(H-\tau-i+\frac12\right)\nabla
    f+\left(2i+\frac52\right)f
\vspace{2mm}\\
\displaystyle \qquad {}=(-1)^{p(k)}\left(2i+k+\frac52\right)f.
\end{array}
\]

\item[5)] Polynomials $f_{k, j}(-H)$ are orthogonal relative to the form
$\langle \cdot, \cdot \rangle_{-i}$; it is a degree $k$ polynomial in $H$ and~$\tau$.

\item[6)] Polynomial $f_{k, 2i}(-H)$ satisfies the difference equation
\[\arraycolsep=0em
\begin{array}{l}
\displaystyle {}\left[(H+1)^2-\left(\tau+\frac12\right)^2\right]\Delta_{2}\nabla_{2} f+4(i+1)(H-i)\nabla_{2}
    f+4i(i+1)f
\vspace{3mm}\\
\displaystyle \qquad {}=\left\{\begin{array}{ll}
(2i+k)(2i+k+2)f&\text{for $k$ even},
\vspace{1mm}\\
(2i+k-1)(2i+k+1)f\quad &\text{for $k$ odd}.
\end{array}\right.
\end{array}
\]
\item[7)] Polynomial $f_{k, 2i+1}(-H)$ satisfies the difference equation
\[\arraycolsep=0em
\begin{array}{l}
\displaystyle    {}\left[(H+1)^2-\left(\tau+\frac12\right)^2\right]\Delta_{2}\nabla_{2}
   f+[(2i+3)(2h-2i-1)-2\tau]\nabla_{2}  f
\vspace{2mm}\\
\displaystyle \qquad {}+(2i+1)(2i+5)f=
\left\{\begin{array}{ll} (2i+k+1)(2i+k+3)f\quad &\text{for $k$ even},
\vspace{1mm}\\
(2i+k)(2i+k+2)f&\text{for $k$ odd}.
\end{array}\right.\end{array}
\]
\item[8)] Polynomial $f_{k, 2i}(H)$ can be expressed via Hahn
polynomials with parameter $\tau$, namely
\[\arraycolsep=0em
\begin{array}{l}
\displaystyle  f_{2k, 2i}(H)= (-1)^k\frac{(i+1)_{k}\left(i+\frac12-\tau\right)_{k}}{k!}
\vspace{3mm}\\
\displaystyle \phantom{ f_{2k, 2i}(H)={}} {}\times
 {}_{3}F_{2}\left(\left.\begin{array}{c}-k, \; k+2i+1, \; \frac12\left(2i+\frac12-\tau-H\right)
\vspace{1mm}\\
 i+1,\;
 i+\frac12-\tau\end{array} \right| \; 1\right),
\vspace{3mm}\\
\displaystyle f_{2k+1, 2i}(H)= \frac12\left[(2k+2i+1)\tau-
\left(\lambda+\frac12\right)^2\right]f_{2k, 2i}(H).
\end{array}
\]
\item[9)] Polynomials $f_{k, 2i+1}(H)$ can be expressed via Hahn
polynomials with parameter $\tau$, namely
\[\arraycolsep=0em
\begin{array}{l}
\displaystyle  f_{2k, 2i+1}(H)= (-1)^k\frac{(i+2)_{k}\left(i+\frac32-\tau\right)_{k}}{k!}
\vspace{3mm}\\
\displaystyle \phantom{f_{2k, 2i+1}(H)={}}{}\times
 {}_{3}F_{2}\left(\left.\begin{array}{c}-k, \; k+2i+2, \; \frac12\left(2i+\frac52-\tau-H\right)
\vspace{1mm}\\
 i+2,\;  i+\frac32-\tau\end{array}\right| \; 1\right),
\vspace{3mm}\\
\displaystyle f_{2k+1, 2i+1}(H)= -\frac{\tau}{2i+k+2}f_{2k, 2i+1}(H).
\end{array}
\]
\end{enumerate}}

{\bf Another theorem of Pinczon.} Recall (Theorem 2.1, B)) that the
superalgebra ${\mathfrak{B}}_{-\frac12}$ is isomorphic to the Weyl algebra $A_1$
considered as {\it super}algebra with generators $P$, $Q$ and relations
$PQ-QP=1$.  The corresponding isomorphism $\theta$ is given by the
formulas
\[
\theta(F)=\frac{1}{\sqrt{2}}P,\qquad  \theta(G)=\frac{1}{\sqrt{2}}Q,\qquad
\theta(H)=\frac{1}{2}(PQ+QP).
\]
As is easy to verify, $\theta(\tau)=0$.

There is a ${\mathbb{Z}}$-grading of $A_1$ such that (having identified $H$
with $\theta(H)$)
\[
A_1={\mathop{\oplus}\limits_{i\in{\mathbb{Z}}}}(A_1)_{i}, \qquad\text{where \ $(A_1)_{i}=\{u\in
A_1\mid [H, u]=iu\}$\  for $i\in{\mathbb{Z}}$}.
\]
For $f, g\in{\mathbb C}[H]$ and $i\geq 0$ set
\begin{equation}
\langle f, g\rangle_{i}= \langle fQ^i, gQ^i\rangle\qquad\text{and}\qquad
\langle f, g\rangle_{-i}= \langle fP^i, gP^i\rangle.\label{2.6.1}
\end{equation}
where $\langle\cdot, \cdot\rangle$ is the bilinear form on $A_1$
defined in Section~4.  Now, for $i\geq 0$ define the polynomials $f_{k, i}$
from the equations
\begin{equation}
\arraycolsep=0em
\begin{array}{l}
f_{k, 2i}Q^{2i}=(\mbox{ad}\; P)^{2k}\left(Q^{2i+2k}\right),
\vspace{2mm}\\
f_{k, 2i+1}Q^{2i+1}=(\mbox{ad}\; P)^{2k+1}\left(Q^{2i+2k+2}\right).
    \end{array}\label{2.6.2}
\end{equation}

Let us endow the algebra ${\mathbb C}[H]$ with a grading by setting $\deg H=1$.

\medskip

\noindent
{\bf Theorem 2.3.} {\it

\vspace{-2mm}

\begin{enumerate}
\topsep0mm
\partopsep0mm
\parsep0mm
\itemsep0mm

\item[1)] $\langle (A_1)_{i} (A_1)_{j}\rangle=0$ for $i\neq j$.

\item[2)] $f_{k, i}(H)$ are polynomials in $H$ and $\tau$ of degree $k$
orthogonal with respect to the form $\langle \cdot, \cdot\rangle_{i}$.

\item[3)] $f_{k, i}(-H)$ are polynomials in $H$ and $\tau$ of degree $k$ orthogonal
with respect to the form $\langle \cdot, \cdot\rangle_{-i}$.

\item[4)] $f_{k, 2i}(H)$ satisfies the difference equation
\[
\left(H+\frac12\right)\Delta f+\left(H-2i-\frac12\right)\nabla f=2kf.
\]

\item[5)] $f_{k, 2i+1}(H)$ satisfies the difference equation
\[
\left(H+\frac12\right)\Delta f+\left(H-2i-\frac32\right)\nabla f=2kf.
\]

\item[6)] The polynomials $f_{k, 2i}(H)$ can be expressed via Meixner
polynomials:
\[
f_{k, 2i}(H)=\frac{(2i+1)_k(2i+k+1)_k}{k!} \cdot
{}_{2}F_{1}\left(\left. \begin{array}{c}-k, \; 2i-H+\frac12 \vspace{1mm}\\
2i+1\end{array}\right|\; 2\right),
\]
where
\[
{}_{2}F_{1}\left(\left. \begin{array}{c} a_{1}, \; a_{2}\vspace{1mm}\\
 b\end{array}\right| \;
z\right)=\mathop{\sum}\limits_{i=0}^\infty
\displaystyle\frac{(a_{1})_{i}(a_{2})_{i} }{(b)_{i}}\,z^i.
\]

\item[7)] The polynomials $f_{k, 2i+1}(H)$ can be expressed via Meixner
polynomials:
\[
f_{k, 2i+1}(H)=\frac{(2i+2)_k(2i+k+2)_{k+1}}{k!} \cdot
{}_{2}F_{1}\left(\left. \begin{array}{c} -k, \; 2i-H+\frac32 \vspace{1mm}\\
2i+2\end{array}\right| \; 2\right).
\]
\end{enumerate}}

{\bf The case of ${\mathfrak{B}}_{\lambda}$ for
$\lambda\in{\mathbb{Z}}_{\geq 0}$} In this case
${\mathfrak{B}}_{\lambda}=\mbox{Mat}(\lambda +1|\lambda )$ and the
image of $\tau$ under the natural homomorphism
$U({\mathfrak{osp}}(1|2))\longrightarrow {\mathfrak{B}}_{\lambda}$ is
a polynomial in $H$.  Therefore, having applied the arguments after
Theorem 2.1 (on the structure of ${\mathfrak{B}}_{\lambda}$) we obtain
an orthogonal basis distinct from the basis of orthogonal polynomials.

To construct orthogonal polynomials, set
\[
U=F, \qquad V=\left(\tau-H+\frac12\right)G.
\]
It is easy to verify that
\begin{equation}
HU-UH=-U, \qquad  HV-VH=V, \qquad  VU-UV=H.\label{2.8.1}
\end{equation}
Relations (\ref{2.8.1}) mean that $U$, $V$, $H$ generate in $U({\mathfrak{osp}}(1|2))$
a subalgebra isomorphic to $U({\mathfrak{sl}}(2))$ considered as a superalgebra
such that $p(U)=p(V)=\overline{1}$, $p(H)=\overline{0}$.  Observe also that the images
of $U$, $V$, and $H$ in ${\mathfrak{B}}_{\lambda}$
generate ${\mathfrak{B}}_{\lambda}$ and are
subject to relations
\begin{equation}
VU=\frac12\left(\lambda(\lambda+1)-H^2+H\right), \qquad
UV=\frac12\left(\lambda(\lambda+1)-H^2-H\right)\label{2.8.2}
\end{equation}
(we have identified $U$, $V$, and $H$ with their images in
${\mathfrak{B}}_{\lambda}$).
The superalgebra ${\mathfrak{B}}_{\lambda}$ is ${\mathbb{Z}}$-graded
${\mathfrak{B}}_{\lambda}=\mathop{\oplus}\limits_{i\in {\mathbb{Z}}}
({\mathfrak{B}}_{\lambda})_{i}$,
where
\[
({\mathfrak{B}}_{\lambda})_{i}=\{u\in {\mathfrak{B}}_{\lambda}\mid [H, u]=iu\}.
\]
Being a matrix superalgebra, ${\mathfrak{B}}_{\lambda}$ possesses an
antiautomorphism, the supertransposition, which in terms of the
generators is given by the formula
\[
H^t=H, \qquad   U^t=-V, \qquad  V^t=U.
\]
The supertrace gives rise to a bilinear form $\langle u, v\rangle
=\mbox{str}\left(uv^t\right)$ on ${\mathfrak{B}}_{\lambda}$.

For $f, g\in{\mathbb C}[H]$ and $i\geq 0$ define the bilinear forms
\begin{equation}
\langle f, g\rangle_{i}= \langle fV^i, gV^i\rangle\qquad\text{and}\qquad
\langle f, g\rangle_{-i}= \langle fU^i, gU^i\rangle.\label{2.8.3}
\end{equation}
Further on, for $i\in{\mathbb{Z}}_{\geq 0}$ set
\begin{equation}\arraycolsep=0em
\begin{array}{l}
f_{0, 2i}V^{2i}=V^{2i}, \quad f_{2, 2i}V^{2i}=\left[U, V^{2i+1}\right], \quad \dots ,
\vspace{2mm}\\
\displaystyle f_{2k, 2i}V^{2i}=\left[U, \frac{1}{H-\left(i+\frac12\right)}\left[U, \left[U,
\frac{1}{H-\left(i+\frac32\right)}\cdots\right.\right.\right.
\vspace{2mm}\\
\left.\left.\left.\displaystyle \phantom{f_{2k, 2i}V^{2i}={}} \frac{1}{H-\left(i+k-\frac32\right)}\left[U, \left[U,
V^{2i+2k-1}\right]\cdots \right]\right]\right]\right].
\end{array}    \label{2.8.4}
\end{equation}
For $u, v\in{\mathfrak{B}}_{\lambda}$ set
\[
\{u, v\}= uv-(-1)^{p(u)(p(v)+1)}vu
\]
and define:
\begin{equation}\arraycolsep=0em
\begin{array}{l}
f_{1, 2i}V^{2i}=\left\{U, V^{2i+1}\right\}, \quad f_{3, 2i}V^{2i}=\left\{U, V^{2i+1}\right\}, \quad \dots ,
\vspace{2mm}\\
\displaystyle f_{2k+1, 2i}V^{2i}=\left\{U, \frac{1}{H-\left(i+\frac12\right)}\left\{U, \left\{U,
\frac{1}{H-\left(i+\frac32\right)}\cdots\right.\right.\right.
\vspace{2mm}\\
\left.\left.\left.\displaystyle \phantom{f_{2k+1, 2i}V^{2i}={}} \frac{1}{H-\left(i+k-\frac32\right)}\left\{U, \left\{U,
V^{2i+2k-1}\right\}\cdots \right\}\right\}\right\}\right\}.
\end{array}    \label{2.8.5}
\end{equation}
Further on, set
\begin{equation}\arraycolsep=0em
\begin{array}{l}
f_{0, 2i+1}V^{2i+1}=V^{2i+1}, \quad f_{2, 2i+1}V^{2i+1}=\left[U, V^{2i+2}\right], \quad \dots ,
\vspace{2mm}\\
\displaystyle
f_{2k, 2i+1}V^{2i+1}=\left[U, \frac{1}{H-\left(i+\frac12\right)}\left[U, \left[U,
\frac{1}{H-\left(i+\frac32\right)}\cdots\right.\right.\right.
\vspace{2mm}\\
\left.\left.\left.\displaystyle \phantom{f_{2k, 2i+1}V^{2i+1}={}} \frac{1}{H-\left(i+k-\frac32\right)}\left[U, \left[U,
V^{2i+2k}\right]\cdots \right]\right]\right]\right],
\end{array}    \label{2.8.6}
\end{equation}
and
\begin{equation}\arraycolsep=0em
\begin{array}{l}
f_{1, 2i+1}V^{2i+1}=\left\{U, V^{2i+2}\right\}, \quad f_{3, 2i+1}V^{2i+1}=\left\{U,
V^{2i+2}\right\}, \quad \dots ,
\vspace{2mm}\\
\displaystyle f_{2k+1, 2i+1}V^{2i+1}=\left(H-\left(i+\frac12\right)\right)f_{2k, 2i+1}.
\end{array}    \label{2.8.7}
\end{equation}

\newpage

\noindent
{\bf Theorem 2.4.} {\it

\vspace{-2mm}

\begin{enumerate}
\topsep0mm
\partopsep0mm
\parsep0mm
\itemsep0mm
\item[1)] $\langle ({\mathfrak{B}}_{\lambda})_{i},
({\mathfrak{B}}_{\lambda})_{j}\rangle=0$ for $i\neq j$.

\item[2)] $f_{l, 2i}$ are orthogonal polynomials of degree $l$ with
respect to the form $\langle\cdot, \cdot \rangle_{2i}$.

\item[3)] Polynomials $f_{2k, 2i}(H)$ satisfy the difference equation
\[
\frac{(H-\lambda)(H+\lambda+1)}{2H-2i+1}\Delta
f+\frac{(H-2i-\lambda-1)(H-2i+\lambda)}{2H-2i-1}\nabla f=2kf.
\]
\item[4)] Polynomials $f_{2k+1, 2i}(H)$ satisfy the difference equation
\[\arraycolsep=0em
\begin{array}{l}
\displaystyle \frac{(H-\lambda)(H+\lambda+1)}{2H-2i+1}\Delta
f+\frac{(H-2i-\lambda-1)(H-2i+\lambda)}{2H-2i-1}\nabla f
\vspace{2mm}\\
\displaystyle \qquad{}+
\frac{\left(\lambda+\frac12\right)^2-i^2}{(H-i)^2-\frac 14}f=(2k+1)f.
\end{array}
\]
\item[5)] Polynomials $f_{l, 2i}(-H)$ are orthogonal with respect to
the form $\langle\cdot, \cdot \rangle_{-2i}$.
\item[6)] Polynomials $f_{l, 2i+1}(H)$ are of degree $l$ and satisfy the
following relations
\[
\langle f_{l, 2i+1}, f_{m, 2i+1}\rangle_{2i+1}\neq 0
\text{ only if }\{l, m\}= \{2k, 2k+1\}_{k\in{\mathbb{Z}}_{\geq 0}}.
\]

\item[7)] Polynomials $f_{l, 2i+1}(-H)$ are of degree $l$ and satisfy
the following relations
\[
\langle f_{l, 2i+1}(-H), f_{m, 2i+1}(-H)\rangle_{-(2i+1)}\neq 0
\text{ only if }\{l, m\}= \{2k, 2k+1\}_{k\in{\mathbb{Z}}_{\geq 0}}.
\]
\item[8)] Polynomials $f_{2k+1, 2i}(-H)$ satisfy the difference equation
\[
\frac{(H-\lambda)(H+\lambda+1)}{2H-2i-1}\Delta
f+\frac{(H-2i-1+\lambda-1)(H-2i-2+\lambda)}{2H-2i-1}\nabla f=2kf.
\]
\item[9)] Polynomials $f_{2k, 2i}$ can be expressed via the dual Hahn
polynomials:
\[
f_{2k, 2i}(H)= (-1)^k\frac{(2i+2)_{k}(\lambda +i+2)_{k}(i-\lambda)_{k}}{k!}
 \cdot
 {}_{3}F_{2}\left(\!\!\!\left.\begin{array}{c}-k, \; i-H, \; H+i\vspace{1mm}\\
\lambda + i+2,\;
 i-\lambda \end{array}\right|\, 1\!\right).
\]
\end{enumerate}}

\setcounter{equation}{0}

\section{Proof: the case $\pbf{\lambda\not\in\left\{-\frac12\right\}\cup{\mathbb{Z}}_{\geq 0}}$}

\noindent
{\bf Lemma 3.1.} {\it Let $A$ be an associative superalgebra generated
by a set $X$.  Denote by $[X, A]$ the set of linear combinations of
the form $\sum [x_{i}, a_{i}]$, where $x_{i}\in X$, $a_{i}\in A$.
Then $[A, A]=[X, A]$.}

\medskip

\noindent
{\bf Proof.} Let us apply the identity \cite[p.~561]{sergeev:Mo}
\begin{equation}
[ab, c]=[a, bc]+\varepsilon(a, bc)[b, ca], \label{3.1.1}
\end{equation}
where $\varepsilon(a, bc)=(-1)^{p(a)(p(b)+p(c))}$.

Namely, let $a=x_1\cdots x_n$; let us perform induction on $n$ to prove that $[a,
A]\subset [X, A]$.  For $n=1$ the
statement is obvious.  If $n>1$, then $a=xa_1$, where $x\in X$ and due
to (\ref{3.1.1}) we have
\[
[a, c]=[xa_{1}, c]=[x, a_{1}c]+\varepsilon(x, a_{1}c)[a_{1}, cx].
\hspace{186pt}\mbox{\rule{3mm}{3mm}}
\]

\noindent
{\bf Lemma 3.2.} {\it Let $A$ be an associative superalgebra and
$a\mapsto a^{t}$ be its antiautomorphism (supertransposition, i.e., it
satisfies $(ab)^{t}=(-1)^{p(a)p(b)}b^{t}a^{t}$ and
$(a^{t})^{t}=(-1)^{p(a)}a$).  Let $L$ be an even invariant functional
on $A$ (like supertrace, i.e., $L([A, A])=0$) such that $L(a^t)=L(a)$
for any $a\in A$.  Define the bilinear form on $A$ by setting
\begin{equation}
\langle u, v\rangle =L\left(uv^t\right) \quad \text{for any} \ \ u, v\in A. \label{3.2.1}
\end{equation}
Then $\langle u, v\rangle =\langle v, u\rangle $ and}
\begin{equation}
\langle [w, u], v\rangle =(-1)^{p(w)(p(u)+1)}\langle u, [w^t, v]\rangle. \label{3.2.2}
\end{equation}

\noindent
{\bf Proof.} Observe first that
\[\arraycolsep=0em
\begin{array}{l}
\langle u, v\rangle =L\left(uv^t\right) =L\left(\left(uv^t\right)^t\right)=
(-1)^{p(u)p(v)}L\left(\left(v^t\right)^tu^t\right)
\vspace{2mm}\\
\phantom{\langle u, v\rangle} {} = (-1)^{p(v)(p(u)+1)}L\left(vu^t\right)
(-1)^{p(v)(p(u)+1)}\langle v, u\rangle .
\end{array}
\]
Since $L$ is even, we see that $\langle u, v\rangle \neq 0$ only if
$p(v)=p(u)$.  But in this case $(-1)^{p(v)(p(u)+1)}$ $=1$.

Further on:
\[\arraycolsep=0em
\begin{array}{l}
\langle [w, u], v\rangle =L\left([w, u]v^t\right)= L\left(\left[w, uv^t\right]-(-1)^{p(w)p(u)}u\left[w,
v^t\right]\right)
\vspace{2mm}\\
\phantom{\langle [w, u], v\rangle} {}= (-1)^{p(w)(p(u)+1)}L\left(u\left[w, v^t\right]\right).
\end{array}
\]
But $[w, v^t]= (-1)^{p(w)+1}\left(\left[w^t, v\right]\right)^t$. Therefore,
\[
\langle [w, u], v\rangle =(-1)^{p(w)(p(u)+1)}L\left(u\left(\left[w^t, v\right]\right)^t\right)=
(-1)^{p(w)(p(u)+1)}\langle u, \left[w^t, v\right]\rangle .\hspace{40.7pt}\mbox{\rule{3mm}{3mm}}
\]

\noindent
{\bf Lemma 3.3.} {\it Set $\{u, v\}=uv-(-1)^{p(u)(p(v)+1)}vu$  and let
$\langle \cdot, \cdot\rangle$ be the bilinear form
as in Lemma~3.2.  Then for $u, v, w\in A$ we have
\[
\langle \{w, u\}, v\rangle =(-1)^{p(w)(p(u)+1)}\langle u, \left\{w^t, v\right\}\rangle .
\]}

\noindent
{\bf Proof.} It is not difficult to verify the following identities:
\begin{equation}
\arraycolsep=0em
\begin{array}{l}
    \{w, u\}v=[w, uv]-(-1)^{p(w)(p(u)+1)}u\{w, v\},
\vspace{2mm}\\
\left\{w, v^t\right\}=-\left(\left\{w^t, v\right\}\right)^t.
\end{array}\label{3.3.1}
\end{equation}
They imply
\[\arraycolsep=0em
\begin{array}{l}
\langle  \{w, u\}, v\rangle =L\left(\{w, u\}v^t\right)=L\left(\left[w,
uv^t\right]-(-1)^{p(u)(p(v)+1)}u\left\{w, v^t\right\}\right)
\vspace{2mm}\\
\phantom{\langle \{w, u\}, v\rangle} {}=-(-1)^{p(u)(p(v)+1)}L\left(u\left\{w, v^t\right\}\right)=
-(-1)^{p(u)(p(v)+1)}L\left(u\left(\left\{w^t, v\right\}\right)^t\right)\vspace{2mm}\\
\phantom{\langle \{w, u\}, v\rangle} {}=-(-1)^{p(u)(p(v)+1)}\langle  u, \left\{w^t, v\right\}\rangle .
\hspace{189.1pt}\mbox{\rule{3mm}{3mm}}
\end{array}
\]

\noindent
{\bf Lemma 3.4.} {\it On ${\mathfrak{B}}_{\lambda}$, there exists a unique, up to
a scalar factor, invariant linear functional, $L$.  It is uniquely
determined by its restriction onto ${\mathbb C}[H]$.  To every functional $L$
on ${\mathbb C}[H]$ assign its generating function
$\varphi_{L}(t)=\sum\limits_{k=0}^\infty\frac{L(H^k)}{k!}t^k$.
 Then for a constant $C\in{\mathbb C}$ we have}
\[
\varphi_{L}(t)=C\frac{e^{(\lambda+1)t}+e^{-\lambda t}}{e^{t}+1}.
\]

\noindent
{\bf Proof.} By Theorem 2.1 the superalgebra ${\mathfrak{B}}_{\lambda}$ is
generated by $G$, $H$, and $F$ subject to relations
\[
\arraycolsep=0em
\begin{array}{l}
{}[H, G]=G,\qquad [H, F]=-F,\qquad  [G, F]=H,\vspace{2mm}\\
\displaystyle   \tau=H+\frac12-2FG,\qquad \tau^2=\left(\lambda+\frac12\right)^2.
\end{array}
\]
Recall that $({\mathfrak{B}}_{\lambda})_{i}=\{u\in {\mathfrak{B}}_{\lambda}\mid [H, u]=iu\}$
for $i\in{\mathbb{Z}}$.
Then for $i\in{\mathbb{Z}}_{\geq 0}$ we have
\[
({\mathfrak{B}}_{\lambda})_{i}=\{fG^i\mid f\in{\mathbb C}[H, \tau]\}, \qquad
({\mathfrak{B}}_{\lambda})_{-i}=\{fF^i\mid f\in{\mathbb C}[H, \tau]\}.
\]
Therefore, $L(({\mathfrak{B}}_{\lambda})_{i})=0$ if $i\neq 0$, so any trace $L$ is
only nonzero on $({\mathfrak{B}}_{\lambda})_{0}={\mathbb C}[H, \tau]$.  To this
restriction assign the generating function
\[
\varphi_{L}(t)=\mathop{\sum}\limits_{k=0}^\infty\frac{L\left(H^k\right)}{k!}t^k
+\frac{\varepsilon}{\lambda+\frac12}\mathop{\sum}\limits_{k=0}^\infty\frac{L\left(\tau
H^k\right)}{k!}t^k
\]
where $\varepsilon^2=1$ and $\varepsilon t=t \varepsilon$.  The following statements are
easy to verify:

i) If $\theta$ is an automorphism of ${\mathbb C}[H, \tau]$ such that
$\theta(H)=H+1$, $\theta(\tau)=-\tau$, and $L_{1}(f)=L(\theta(f))$,
then $\varphi_{L_{1}}=e^t\overline{\varphi_{L}}$, where
$\overline{\varphi_{L}}= \varphi_{L, \overline{0}}-\varepsilon \varphi_{L, \overline{1}}$ for
each $\varphi_{L}= \varphi_{L, \overline{0}}+\varepsilon \varphi_{L, \overline{1}}$, where
$\varphi_{L, \overline{0}}$ and $\varphi_{L, \overline{1}}$ are formal power series in~$t$.

ii) If $L_{2}(f)=L(\tau f)$, where $f\in {\mathbb C}[H, \tau]$, then
\[
\varphi_{L_{2}}=\left(\lambda +\frac12\right)\varepsilon\varphi_{L}.
\]

iii) If $P\in {\mathbb C}[H]$, and $L_{3}(f)=L(Pf)$ for $f\in {\mathbb C}[H, \tau]$,
then
\[
\varphi_{L_{3}}=P\left(\frac{d}{dt}\right)\varphi_{L}.
\]

Making use of these statements, let us calculate the generating
function for the restriction of the functional $L$ onto ${\mathbb C}[H,
\tau]$. For $f\in {\mathbb C}[H, \tau]$ we have
\[\arraycolsep=0em
\begin{array}{l}
\displaystyle  {}[F, fG]=FfG+fGF=\theta(f)FG+fGF=\theta(f)\frac12\left(H-\tau+\frac12\right)
\vspace{3mm}\\
\displaystyle \phantom{[F, fG]={}} {}+\frac12f\left(H-\tau+\frac12\right)=
\theta\left(\frac12f\left(H-\tau+\frac12\right)\right)+\frac12f\left(H-\tau+\frac12\right).
\end{array}
\]
Hence, $L\left(\theta\left(\frac12f\left(H-\tau+\frac12\right)\right)+\frac12f\left(H-\tau+\frac12\right)
\right)=0$.
Therefore, thanks to i)--iii) we have
\begin{equation}
\left (\frac{d}{dt}-\frac12+\left(\lambda+\frac12\right)\varepsilon\right)
\left(e^t\bar\varphi_{L}+
\varphi_{L}\right)=0.
\label{3.4.1}
\end{equation}
For $\varphi_{L}= \varphi_{L, \overline{0}}+\varepsilon \varphi_{L, \overline{1}}$ we obtain a
system:
\[\arraycolsep=0em
\begin{array}{l}
\displaystyle \left (\frac{d}{dt}-\frac12\right)\left(e^t+1\right)\varphi_{L,
\overline{0}}+\left (\lambda+\frac12\right)\left(e^t-1\right)
\varphi_{L, \overline{1}}=0
\vspace{3mm}\\
\displaystyle \left (\frac{d}{dt}-\frac12\right)\left(e^t-1\right) \varphi_{L, \overline{1}}+
\left (\lambda+\frac12\right)\left(e^t+1\right)\varphi_{L, \overline{0}}=0.
\end{array}
\]
Hence,
\[
\left
(\frac{d}{dt}-\frac12\right)^2\left(e^t-1\right) \varphi_{L, \overline{1}}+
\left (\lambda+\frac12\right)^2\left(e^t-1\right)\varphi_{L, \overline{1}}=0
\]
or, even simpler,
\[
\left(\frac{d}{dt}-\lambda-1\right)
\left(\frac{d}{dt}+\lambda\right)\left(e^t-1\right) \varphi_{L, \overline{1}}=0.
\]
This implies that
\[
\varphi_{L, \overline{1}}=c\frac{e^{(\lambda+1)t}-e^{-\lambda t}}{e^t-1}\qquad
\text{and} \qquad\varphi_{L, \overline{0}}=c\frac{e^{(\lambda+1)t}+e^{-\lambda
t}}{e^t+1}.
\]
Since $\varphi_{L, \overline{1}}$ is uniquely recovered from $\varphi_{L,
\overline{0}}$, we see that $L$ is uniquely recovered by its restriction
onto ${\mathbb C}[H]$. This proves uniqueness.

Let us prove existence of $L$. It suffices to prove that $1\not
\in [{\mathfrak{B}}_{\lambda}, {\mathfrak{B}}_{\lambda}]$. Indeed, by Lemma~3.1
\[
[{\mathfrak{B}}_{\lambda}, {\mathfrak{B}}_{\lambda}]=
[F, {\mathfrak{B}}_{\lambda}]+[G, {\mathfrak{B}}_{\lambda}].
\]
Hence,
\[\arraycolsep=0em
\begin{array}{l}
  {}[{\mathfrak{B}}_{\lambda}, {\mathfrak{B}}_{\lambda}]\cap{\mathbb C}[H, \tau]
=([F, {\mathfrak{B}}_{\lambda}]+[G,  {\mathfrak{B}}_{\lambda}])\cap{\mathbb C}[H, \tau]
\vspace{2mm}\\
\displaystyle \qquad {}= [F, ({\mathfrak{B}}_{\lambda})_{1}]\cap{\mathbb C}[H,
\tau]=\mbox{Span}\left(\theta\left(f\left(H+\tau-\frac12\right)\right)+
f\left(H+\tau-\frac12\right)\right)
\end{array}
\]
for any $f\in{\mathbb C}[H, \tau]$.

But $\deg [F, fG]\geq 1$ (recall that $\deg H=2$, $\deg \tau =1$),
so $1\not
\in [{\mathfrak{B}}_{\lambda}, {\mathfrak{B}}_{\lambda}]$. \hfill\rule{3mm}{3mm}

\medskip

\noindent
{\bf Lemma 3.5.} {\it Let $L$ be a linear functional on
${\mathfrak{B}}_{\lambda}$ determined in Lemma 3.4 and normed
so that
$L(1)=1$. Define an antiautomorphism of ${\mathfrak{B}}_{\lambda}$ by setting
\[
H^t=H, \qquad G^t=-F, \qquad F^t=G.
\]
Then
\vspace{-2mm}

\begin{enumerate}
\topsep0mm
\partopsep0mm
\parsep0mm
\itemsep0mm

\item[i)] $\left(u^t\right)^t=(-1)^{p(u)}u$ for any $u\in
{\mathfrak{B}}_{\lambda}$;

\item[ii)] $L\left(u^t\right)=L(u)$ for any $u\in
{\mathfrak{B}}_{\lambda}$.
\end{enumerate}}

\noindent
{\bf Proof.} i) Induction on $n$, where $u=x_{1}\dots x_{n}$ and
$x_{i}\in\mbox{Span}\,(H, G, F)$ for each $i$.

For $n=1$ the statement is obvious. Let $u=u_{1}u_{2}$ and let for
$u_{1}$ and for $u_{2}$ the statement be true. Then
\[\arraycolsep=0em
\begin{array}{l}
\left(u^t\right)^t=\left((u_{1}u_2)^t\right)^t=\left((-1)^{p(u_1)p(u_2)}u_2^tu_1^t\right)^t
\vspace{2mm}\\
\phantom{\left(u^t\right)^t} {}=
\left(u_1^t\right)^t\left(u_2^t\right)^t=(-1)^{p(u_1)+p(u_2)}u_1u_2=(-1)^{p(u)}u.
\end{array}
\]

ii) Let us represent $u\in {\mathfrak{B}}_{\lambda}$ in the form $u=u_{0}+L(u)$,
where $u_{0}\in \mbox{Ker}\; L$. But we know that $\mbox{Ker}\; L=[{\mathfrak{B}}_{\lambda},
{\mathfrak{B}}_{\lambda}]$. So
\[
(\mbox{Ker}\; L)^t=([{\mathfrak{B}}_{\lambda},
{\mathfrak{B}}_{\lambda}])^t= [{\mathfrak{B}}_{\lambda},
{\mathfrak{B}}_{\lambda}]= \mbox{Ker}\; L.
\]
Hence, $u_{0}^t\in \mbox{Ker}\; L$ and $u^t=u_{0}^t+L(u)$. Thus, $L\left(u^t\right)=L(u)$.
\hfill \rule{3mm}{3mm}

\medskip

Therefore, we can define a bilinear form $\langle u, v\rangle=L\left(u^tv\right)$
on ${\mathfrak{B}}_{\lambda}$.  By Lemma~3.2 we have: $\langle u,
v\rangle=(-1)^{p(u)p(v)}\langle v, u\rangle$ and
\[
\langle [w, u], v\rangle=(-1)^{(p(u)+1)p(w)}\langle u, \left[w^t, v\right]\rangle.
\]

{\bf Proof of heading 1 of Theorem 2.2.} Let
$u\in({\mathfrak{B}}_{\lambda})_{i}$, $u\in({\mathfrak{B}}_{\lambda})_{j}$.  Then
\[
i\langle u, v\rangle=\langle [H, u], v\rangle=\langle u, [H,
v]\rangle=j\langle u, v\rangle.
\]
Therefore, if $i\neq j$, then $\langle u, v\rangle=0$.

{\bf Proof of heading 2 of Theorem 2.2.} By~\cite{sergeev:Pi} (see also
\cite{sergeev:GL}), there is an expansion
${\mathfrak{B}}_{\lambda}=\mathop{\oplus}\limits_{n\in{\mathbb{Z}}_{\geq 0}}
(\Pi)^{n}\left({\mathcal{L}}^n\right)$,
where ${\mathcal{L}}^n$ is an irreducible highest weight module over
${\mathfrak{osp}}(1|2)$ with even highest weight vector, $\Pi$ is the change of
parity functor and where ${\mathfrak{B}}_{\lambda}$ is considered as
${\mathfrak{osp}}(1|2)$-module with respect to the adjoint representation.  It is
easy to verify that ${\mathcal{L}}^{2n}$ is generated by the highest weight
vector $G^{2n}$, whereas $\Pi\left({\mathcal{L}}^{2n+1}\right)$ is generated by
$G^{2n+1}\tau$.  Hence, $f_{k, 2i}G^{2i}\in {\mathcal{L}}^{2i+k}$.

Making use of Lemma 3.5, it is not difficult to verify that
\[
\langle \Omega *u, v\rangle= \langle u, \Omega *v\rangle,
\]
where $\Omega$ is defined in Lemma 2.1 and $*$ denotes the adjoint
action.  This immediately implies that $\langle {\mathcal{L}}^p, {\mathcal{L}}^q\rangle=0$
if $p\neq q$ and, therefore, $\langle f_{k, 2i}, f_{l, 2i}\rangle_{2i}=0$ if
$k\neq l$.

Let us show now that $\deg f_{k, 2i}=k$. (Recall again that $\deg H=2$,
$\deg\tau=1$.) It is easy to verify that
\[\arraycolsep=0em
\begin{array}{l}
\displaystyle {}[F, [F, fG^p]]=[fG^p,  Y]=\frac14\Biggl\{\left((H+1)^2-\left(\tau+\frac12\right)^2\right)\Delta_{2}f
\vspace{2mm}\\
\displaystyle \phantom{[F, [F, fG^p]]=[fG^p,  Y]={}}{}+
\left(p(2H-p+2)-\tau+(-1)^p\tau\right)\Biggr\}G^{p-2}=\tilde
fG^{p-2}.
\end{array}
\]
This implies $\deg f=\deg \tilde f+2$.  If $k=0$, then $f_{0,
2i}G^{2i}=G^{2i}$ and $f_{0, 2i}=1$.  Formulas (\ref{2.4.2}) imply that
$f_{0, 2i}G^{2i}=G^{2i}=\left[F, \left[F, f_{k, 2i+2}G^{2i+2}\right]\right]$, hence, $\deg
f_{k+2, 2i}= \deg f_{k, 2i+2}+2$.

We similarly prove that $\langle f_{k, 2i+1}, f_{l,
2i+1}\rangle_{2i+1}=0$ if $k\neq l$ and $\deg f_{k, 2i+1}=k$.

{\bf Proof of heading 3 of Theorem 2.2.} For any $f\in{\mathbb C}[H,
\tau]$ set
\[
\arraycolsep=0em
\begin{array}{l}
\Delta f(H)=f(H+1)-f(H), \qquad \nabla f(H)=f(H)-f(H-1),
\vspace{1mm}\\
\Delta \tau=-2\tau, \qquad \nabla \tau=2\tau.
\end{array}
\]
The following identities are easy to check
\begin{equation}
{}[F, f]=\Delta f\cdot F=F\nabla f,\qquad [G, f]=-G\Delta f=-\nabla fG.\label{3.8.1}
\end{equation}
Moreover,
\begin{equation}
\arraycolsep=0em
\begin{array}{l}
{}\left[F, G^{2i}\right]=iG^{2i-1},\qquad \left[F, G^{2i+1}\right]=(H-i)G^{2i},
\vspace{2mm}\\
\displaystyle FG=\frac12\left(H-\tau+\frac12\right), \qquad
GF=\frac12\left(H+\tau-\frac12\right),
\vspace{2mm}\\
Ff(H)=f(H+1)F,\qquad Gf(H)=f(H-1)G, \vspace{2mm}\\
 F\tau=-\tau F,\qquad G\tau=-\tau FG.
\end{array}\label{3.8.2}
\end{equation}
Let us calculate the results of the adjoint action of $\tau$ on
$fG^{2i}$, where $f\in{\mathbb C}[H, \tau]$.  From the explicit
expression of $\tau$ (Lemma~2.2) we deduce
\[\arraycolsep=0em
\begin{array}{l}
\displaystyle     \tau*\left(fG^{2i}\right)=\frac12fG^{2i}+\left[H, fG^{2i}\right]-2\left[F\left[G, fG^{2i}
\right]\right]
\vspace{2mm}\\
\displaystyle \phantom{ \tau*\left(fG^{2i}\right)} {}= \left(2i+\frac12\right)fG^{2i}+2\left[F,\nabla f\cdot
G^{2i+1}\right]
\vspace{2mm}\\
\displaystyle \phantom{ \tau*\left(fG^{2i}\right)} {} =
2\Delta\nabla f\cdot FG\cdot G^{2i}+2\nabla f\cdot\left[F, G^{2i+1}\right]+
\left(2i+\frac12\right)fG^{2i}
\vspace{2mm}\\
\displaystyle \phantom{ \tau*\left(fG^{2i}\right)} {}
= \left[\left(H-\tau+\frac12\right)\Delta\nabla f+2(H-i)\Delta f+\left(2i+\frac12\right)f\right]G^{2i}.
\end{array}
\]

On the other hand, if $fG^{2i}\in{\mathcal{L}}^{2i+k}$, then $\tau*\left(fG^{2i}\right)=
c\cdot fG^{2i}$ because in any irreducible ${\mathfrak{osp}}(1|2)$-module $\tau$
acts as a scalar multiple of the parity operator $P$, i.e., an
operator such that $P(v)=(-1)^{p(v)}v$ for any $v\in V$.

Operator $\tau$ acts on the highest weight vector of ${\mathcal{L}}^{2i+k}$ as
multiplication by $2i+k+\frac12$.  Observe that $p\left(f_{k,
2i}G^{2i}\right)=p\left(G^{2i+k}\right)$ if $k$ is even and $p\left(f_{k,
2i}G^{2i}\right)=p\left(G^{2i+k}\tau\right)+\overline{1}$ if $k$ is odd.

Therefore,
\[
\tau* \left(f_{k, 2i}G^{2i}\right)=(-1)^k\left(2i+k+\frac12\right)f_{k,
2i}G^{2i}.
\]

{\bf Proof of heading 4 of Theorem 2.2.} Let us calculate the
results of the adjoint action of $\tau$ on $fG^{2i+1}$, where
$f\in{\mathbb C}[H, \tau]$. We obtain
\[
\arraycolsep=0em
\begin{array}{l}
\displaystyle   \tau* \left(fG^{2i}\right)=\left(2i+\frac32\right)fG^{2i}+2\left[F, \nabla
(fG^{2i+2})-2fG^{2i+2}\right]
\vspace{2mm}\\
\displaystyle \qquad {}=
2\Delta\nabla f\cdot FG\cdot G^{2i+1}+2\nabla f\cdot \left[F,
G^{2i+2}\right]-4\nabla f\cdot FG\cdot G^{2i+1}-4f\left[F, G^{2i+2}\right]
\end{array}
\]
\[
\arraycolsep=0em
\begin{array}{l}
\displaystyle \qquad {}+\left(2i+\frac32\right)fG^{2i+1}=
\left(H-\tau+\frac12\right)\Delta\nabla f\cdot G^{2i+1}+2(i+1)\nabla f\cdot G^{2i+1}
\vspace{2mm}\\
\displaystyle \qquad {}-2\left(H-\tau+\frac12\right)
(\Delta\nabla f+\nabla
f)G^{2i+1}-4(i+1)fG^{2i+1}+\left(2i+\frac32\right)fG^{2i+1}
\vspace{2mm}\\
\displaystyle \qquad {}=
-\left\{\left(H-\tau+\frac12\right)\Delta\nabla f+2\left(H-\tau-i-\frac12\right)
\nabla f+\left(2i+\frac52\right)f\right\}G^{2i+1}.
\end{array}
\]

On the other hand, as in the proof of heading 3, we see that
\[
\tau*\left(f_{k, 2i+1}G^{2i+1}\right)= c\cdot f_{k, 2i+1}G^{2i+1}.
\]
But $f_{k,
2i+1}G^{2i+1}\in\left\{\!\!\begin{array}{ll}{\mathcal{L}}^{2i+k+2}&\text{
for $k$ even}\\ {\mathcal{L}}^{2i+k}&\text{ for $k$
odd}\end{array}\right.$ and the parity of $f_{k, 2i+1}G^{2i+1}$
coincides with that of the highest weight vector of
${\mathcal{L}}^{2i+k+2}$ if $k$ is even, and is opposite if $k$ is
odd; so
\[
\tau* \left(f_{k,
2i+1}G^{2i+1}\right)=-(-1)^{p(k)}\left(2i+k+\frac52\right)f_{k,
2i+1}G^{2i+1}.
\]

{\bf Proof of heading 5 of Theorem 2.2.} Let $\theta$ be an
automorphism of ${\mathfrak{B}}_{\lambda}$ given on generators as
follows:
\begin{equation}
 \theta(G)=\sqrt{-1}F, \qquad \theta(F)=\sqrt{-1}G, \qquad\theta(H)=-H.\label{3.10.1}
\end{equation}
Let $L$ be a functional on ${\mathfrak{B}}_{\lambda}$ defined in Lemma
3.4.  Since $L$ is unique, up to a scalar factor, invariant linear
functional on ${\mathfrak{B}}_{\lambda}$, it follows that
$\mbox{Ker}\; L=[{\mathfrak{B}}_{\lambda}, {\mathfrak{B}}_{\lambda}]$. 
Hence, $\theta(\mbox{Ker}\; L)=\mbox{Ker}\; L$ and
$L(\theta(u))=L(u)$, if $u\in {\mathfrak{B}}_{\lambda}$.  Therefore,
\[\arraycolsep=0em
\begin{array}{l}
    \langle f, g\rangle_{-i}=\langle fF^i,
gF^i\rangle=L\left(fF^i\left(gF^i\right)^t\right)=(-1)^{i(i-1)/2}L\left(fgF^iG^i\right)
\vspace{2mm}\\
\phantom{\langle f, g\rangle_{-i}} {} =
(-1)^{i(i-1)/2}L\left(\theta\left(fgF^iG^i\right)\right)=
(-1)^{i(i-1)/2}L\left(\theta(f)G^iF^i\theta(g)\right)\left(\sqrt{-1}\right)^{2i}
\vspace{2mm}\\
\phantom{\langle f, g\rangle_{-i}} {}= (-1)^iL\left(\theta(f)G^i\right)\left(\theta(g)G^i\right)^t=(-1)^i\langle\theta(f),
\theta(g)\rangle_{-i}.
\end{array}
\]
But $\theta(H)=-H$ and
\[
\theta(\tau)=\theta\left(H+\frac12-2FG\right)=-H+\frac12+2GF=\tau.
\]

{\bf Proof of heading 6 and 7 of Theorem 2.2.} Recall that
$\Delta_{2}(H)=\nabla_{2}(H)=H$ and $\Delta_{2}(\tau)=0$. Moreover,
$X=G^2$, $Y=F^2$, and $H$ span a Lie algebra isomorphic to ${\mathfrak{sl}}(2)$
and the following relations hold
\begin{equation}
XY=\frac14\left(\left(\tau+\frac12\right)^2-(H-1)^2\right), \qquad
YX=\frac14\left(\left(\tau+\frac12\right)^2-(H+1)^2\right).\label{3.11.1}
\end{equation}
It is easy to verify that for $f\in{\mathbb C}[H, \tau]$ we have
\[\arraycolsep=0em
\begin{array}{l}
{}[X, f]=-X\Delta_{2}f=-\nabla_2fX,\qquad [Y, f]=\Delta_{2}fY=Y\nabla_2f,
\vspace{2mm}\\
\displaystyle {}\left[Y, G^{2i}\right]\!=-i(H-i+1)G^{2i-2},\quad \left[Y,
G^{2i+1}\right]\!=\frac12\left(\tau-(2i+1)\!\left(H-i+\frac12\right)\right )G^{2i-1}\!.
\end{array}
\]
Let us compute the result of the adjoint action of the Casimir
operator $\omega=H^2+2H+4XY$ from
$U({\mathfrak{sl}}(2))\subset {\mathfrak{B}}_{\lambda}$ on
$fG^{2i}$.  We have
\[\arraycolsep=0em
\begin{array}{l}
    \omega*\left(fG^{2i}\right)=\left[H, \left[H, fG^{2i}\right]\right]
+2\left[H, fG^{2i}\right]+4\left[Y, \left[X, fG^{2i}\right]\right]
\vspace{2mm}\\
\qquad {}= (4i^2+4i)fG^{2i}-4\left[Y,  \Delta_2fG^{2i+2}\right]
\vspace{2mm}\\
\qquad {}= 4i(i+1)fG^{2i}-4\Delta_2\nabla_2fYXG^{2i}+
4i(H-i)\Delta_2fG^{2i}
\vspace{2mm}\\
\displaystyle \qquad {}=
\!\left[(H+1)^2\!-\!\left(\tau+\frac12\right)^2\right]\!\Delta_2\nabla_2fG^{2i}
+4(i+1)(H\!-\!i)\nabla_2fG^{2i}+
4i(i+1)fG^{2i}.
\end{array}
\]
As ${\mathfrak{sl}}(2)$-module, ${\mathcal{L}}^{2i+k}=L^{2i+k}\oplus \Pi\left(L^{2i+k-1}\right)$,
where $L^{m}$ is the irreducible (finite dimensional) ${\mathfrak{sl}}(2)$-module
with highest weight $m$.  As is easy to calculate, $\omega$ acts on
$L^{m}$ as multiplication by $m(m+2)$.

We have
\[
f_{k, 2i}\in\left\{\begin{array}{ll} L^{2i+k}\subset {\mathcal{L}}^{2i+k}&\text{for
$k$ even}, \vspace{1mm}\\
L^{2i+k-1}\subset {\mathcal{L}}^{2i+k}\quad &\text{for $k$ odd.}\end{array}\right.
\]

Now, let us compute the action of $\omega$ on $fG^{2i+1}$:
\[\arraycolsep=0em
\begin{array}{l}
\displaystyle     \omega*\left(fG^{2i+1}\right)=    \Biggl\{
\left[(H+1)^2-\left(\tau+\frac12\right)^2\right]\Delta_2\nabla_2f
\vspace{2mm}\\
\phantom{\omega*\left(fG^{2i+1}\right)={}} {}+ \left[(2i+3)(2H-2i-1)-2\tau\right]\nabla_2f+
(2i+1)(2i+5)f\Biggr\}G^{2i+1}.
\end{array}
\]

{\bf Proof of heading 8 and 9 of Theorem 2.2.} Since
$\Delta_{2}\nabla_{2}=\Delta_{2}-\nabla_{2}$, we can express the left
hand side of the equation of heading 6) as
\[
\left [(H+1)^2-\left(\tau+\frac12\right)^2\right ]\Delta_2f+
\left [(2H-2i-1)^2-\left(\tau+\frac12\right)^2\right ]\nabla_2f.
\]
Making the change $x=\frac12H+\frac 12\tau-\frac14-i$ we reduce the
above equation to the form
\[\arraycolsep=0em
\begin{array}{l}
(x-N)(x+\alpha+1)\Delta \varphi-x(x-\beta-N-1)\nabla\varphi
\vspace{2mm}\\
\qquad \qquad = \left\{ \begin{array}{ll}
\displaystyle \frac{k}{2}\left(\frac{k}{2}+\alpha+\beta+1\right)\varphi&\text{for $k$ even},
\vspace{2mm}\\
\displaystyle \frac{k-1}{2}\left(\frac{k-1}{2}+\alpha+\beta+1\right)\varphi\quad &\text{for $k$
odd,}\end{array}\right.
\end{array}
\]
where $N=\tau-\frac12-i$, $\alpha=\beta=i$ and
$\varphi(x)=f\left(2x+i+\frac12-\tau\right)$.

Thanks to \cite[p.~30]{sergeev:KS}, we know that one of the solutions of the
above equation is equal to ${}_{3}F_{2}\left(\left.\begin{array}{c}-l,\;
l+\alpha+\beta+1, \; -x\\ \alpha+1, \; -N\end{array}\right|\, 1\right)$.
(In~\cite{sergeev:KS} it is supposed that $N$ is a positive integer and $x$ is
real, but one can clearly assume that $N$ and $x$ belong to any
commutative ring.)  Thus,
\[
\varphi(x)=c\cdot {}_{3}F_{2}\left(\left.\begin{array}{c}-l,\;
l+\alpha+\beta+1, \; -x\\ \alpha+1, \; -N\end{array}\right| \,
1\right) \quad\text{for some}\ \ c\in{\mathbb C}[\tau].
\]
To calculate the exact value of the constant $c$, it suffices to
compute the leading coefficient of the polynomial $f_{k, 2i}$.
Formula (\ref{2.4.2}) implies that
\begin{equation}\arraycolsep=0em
\begin{array}{l}
    f_{k, 2i}G^{2i}=\left[F, \left[F, f_{k-2, 2i+2}G^{2i+2}\right]\right]=\left[f_{k-2,
2i+2}G^{2i+2}, Y\right]
\vspace{2mm}\\
\displaystyle  {}=
\frac14\left\{ \!
\left[(H+1)^2\!-\!\left(\tau+\frac12\right)^2\right]\!\Delta_2f_{k-2, 2i+2}+
2i(2H-2i+2)f_{k-2, 2i+2}\right\}G^{2i+2}.
\end{array}\label{3.12.1}\hspace{-5mm}
\end{equation}
In particular,
\[
f_{k-2, 2i+2}=1, \qquad f_{1, 2i+2}G^{2i}=\left[F,
G^{2i+1}\tau\right]=\left (\frac{2i+1}{2}\tau-
\frac{\left(\lambda+\frac12\right)^2}{2}\right )G^{2i}.
\]

Let $f_{k, 2i}=a_{k, 2i}H^{[k/2]}+\cdots$, where $[x]$ denotes the
integer part of~$x$.  Then (\ref{3.12.1}) implies that
\begin{equation}
a_{k, 2i}=\left (\frac12\left [\frac{k-2}{2}\right ]+i+1\right )a_{k-2, 2i+2}. \label{3.12.2}
\end{equation}
Formula (\ref{3.12.2}) implies
that
\[
a_{2l, 2i}=\frac{(2i+l+1)_{l}}{2^l}, \qquad a_{2l+1,
2i}=\left[(2i+l+1)\tau-\left(\lambda+\frac12\right)^2\right]\times \frac{(2i+l+1)_{l}}{2^{l+1}}.
\]
Since the coefficient of the leading power of $x$ in
${}_{3}F_{2}\left(\left.\begin{array}{c}-l, \; l+\alpha+\beta+1, \; -x\\ \alpha+1,\;
-N\end{array} \right| \,  1\right)$ is equal to
\[
\frac{(-1)^ll!(l+\alpha+\beta+1)_{l}}{(\alpha +1)_{l}(-N)_{l}},
\]
we  deduce that
\[
c=\frac{(-1)^l(\alpha +1)_{l}(-N)_{l}}{l!(l+\alpha+\beta+1)_{l}}\cdot
\frac{(2i+l+1)_{l}(-N)_{l}}{2^l}2^l
\]
which leads to formulas of heading 8.

Similar calculations show that
\[
a_{2l, 2i+1}=\frac{(2i+l+2)_{l}}{2^l}, \qquad a_{2l+1,
2i+1}=- \frac{(2i+l+1)_{l}}{2^{l}}\tau.
\]
This leads to formulas of heading 9.

\setcounter{equation}{0}

\section{Proof for $\pbf{\lambda =-\frac12}$}

We will stick to notations introduced after Theorem 2.2 concerning
``another theorem of Pinczon''.

\medskip

\noindent
{\bf Lemma 4.1.} {\it On Weyl algebra $A_{1}={\mathbb C}[P, Q]$, there exists
a unique up to a constant factor invariant linear functional $L$.  It
is uniquely determined by its restriction onto ${\mathbb C}[H]$.  The
generating function of $L$ is of the form}
\[
\varphi_{L}(t)=c\frac{e^{t/2}}{1+e^t}\qquad\text{for}\ \ c\in{\mathbb C}.
\]

\smallskip

\noindent
{\bf Proof.} Since $(A_{1})_{i}=\{u\in A_{1}\mid [H, u]=iu\}$ for
$i\in{\mathbb{Z}}$, we see that $L((A_{1})_{i})=0$ for $i\neq 0$.  So $L$ is
uniquely determined by its restrictions onto $(A_{1})_{0}={\mathbb C}[H]$.

Further,
\[
{}[A_{1}, A_{1}]\cap {\mathbb C}[H]=[P, (A_{1})_{1}]\cap {\mathbb C}[H]=\mbox{Span}([P, fQ])
\]
where $f\in{\mathbb C}[H]$. But
\[\arraycolsep=0em
\begin{array}{l}
[P, fQ]=PfQ+fPQ=f(H+1)PQ+f(H)QP
\vspace{2mm}\\
\displaystyle \phantom{[P, fQ]} {}=\left(H+\frac12\right)f(H+1)+\left(H-\frac12\right)f(H).
\end{array}
\]
And, since $\deg [P, fQ]\geq 1$ (here we assume that $\deg H=1$), it
follows that $1\not \in [A_{1}, A_{1}]$ which proves the existence of
$L$. Further on,
\[
L\left(\left(H+\frac12\right)f(H+1)+\left(H-\frac12\right)f(H)\right)=0
\]
wherefrom, as in Lemma 3.4, we deduce that
\[
\left(\frac{d}{dt}-\frac12\right)\left(1+e^t\right)\varphi_{L}=0
\]
and the desired form of $\varphi_{L}$. \hfill \rule{3mm}{3mm}

\medskip

\noindent
{\bf Lemma 4.2.} {\it Define an automorphism of $A_{1}={\mathbb C}[P, Q]$ by
setting
\[
H^t=h, \qquad Q^t=-P, \qquad P^t=Q.
\]
Then

~i) $\left(u^t\right)^t=(-1)^{p(u)}u$;

ii) $L(u^t)=L(u)$, where $u\in A_{1}$.}

\medskip

\noindent
Proof is similar to that of Lemma 3.5.

Thus, the form $\langle\cdot, \cdot\rangle$ is supersymmetric and
invariant.

{\bf Proof of Theorem 2.3.} Heading 1 is proved as heading~1 of
Theorem 2.2.  Heading~2 is proved as heading~2 of Theorem 2.2 with the
help of decomposition $A_1
=\mathop{\oplus}\limits_{n\in{\mathbb{Z}}_{\geq
0}}{\mathcal{L}}^{2n}$, where ${\mathcal{L}}^{2n}$ is an irreducible
highest weight module over ${\mathfrak{osp}}(1|2)$ with even highest
weight vector.  Heading~3 is proved as heading~5 of Theorem 2.2 with
the help of automorphism~(\ref{3.10.1}), where $F=P$ and $G=Q$.  The
difference equations for $f_{k, j}$ follows from the study of the
result of application of $\tau=H+\frac12-2FG$ to $f_{k, j}G^j$ under
the adjoint action of ${\mathfrak{osp}}(1|2)$ on $A_{1}$ and arguments
similar to those from the proof of heading~3 of Theorem 2.2. 
Statements of headings~6 and~7 are results of comparison of difference
equations in headings 4 and 5 with corresponding equations
in~\cite{sergeev:KS} and calculation of the leading terms.

\setcounter{equation}{0}

\section{Proof for the case $\pbf{\lambda\in {\mathbb{Z}}_{\geq 0}}$}

As was observed in proof after formula (2.13), the elements $U=F$ and
$V=\left(\tau-H+\frac12\right)G$ generate in $U({\mathfrak{osp}}(1|2))$ a subalgebra
isomorphic to $U({\mathfrak{sl}}(2))$, considered as a superalgebra with
nontrivial odd part.

It is also convenient to consider $U({\mathfrak{sl}}(2))$ per se, not as a
subalgebra of $U({\mathfrak{osp}}(1|2))$. We mean the following.

Let ${\mathfrak{sl}}(2)=\mbox{Span}(X, H, Y)$ with relations (\ref{1.1.1}).  Consider
$U({\mathfrak{sl}}(2))$ as a superalgebra with parity given by the formula
$p(X)=p(Y)=\overline{1}$ (hence, $p(H)=\overline{0}$). Set
\[
U=\frac{1}{\sqrt{2}}Y, \qquad V=\frac{1}{\sqrt{2}}X, \qquad H\mapsto
\frac12 H,
\]
we, clearly, have
\begin{equation}
HU-UH=-U, \qquad HV-VH=V, \qquad VU-UV=H.\label{5.1.1}
\end{equation}
The Casimir operator, being even, remains the same:
\[
\Omega=H^2-H+2VU=H^2+H+2UV.
\]
Therefore,
\begin{equation}
\arraycolsep=0em
\begin{array}{l}
\displaystyle VU=\frac12\left(\Omega-\left(H^2-H\right)\right), \qquad
UV=\frac12\left(\Omega-\left(H^2+H\right)\right),
\vspace{2mm}\\
{}[V, U]=\Omega-H^2.
\end{array}\label{5.1.2}
\end{equation}
(Recall that $[U, V]=UV+VU$.)

Let ${\mathfrak{A}}_{\lambda}$ be a quotient of $U({\mathfrak{sl}}(2))$, as in (\ref{1.1.2}). Then
${\mathfrak{B}}_{\lambda}\simeq {\mathfrak{A}}_{2\lambda}$ for any
$\lambda\in{\mathbb{Z}}_{\geq 0}$.

Set ${\mathfrak{C}}_{\lambda}:= {\mathfrak{A}}_{2\lambda}$
for any $\lambda\in{\mathbb C}$; in other words, ${\mathfrak{C}}_{\lambda}$ is generated
by odd indeterminates $U$ and $V$ subject to relations
\begin{equation}
\arraycolsep=0em
\begin{array}{l}
    HU-UH=-U, \qquad HV-VH=V, \qquad VU-UV=H,
\vspace{2mm}\\
\displaystyle VU=\frac12\left(\lambda(\lambda+1)-\left(H^2-H\right)\right), \qquad
UV=\frac12\left(\lambda(\lambda+1)-\left(H^2+H\right)\right),
\end{array}\label{5.1.3}
\end{equation}
and one more relation:
\[
V^{2|\lambda|}=0 \quad \text{for$ \ \ \lambda\in\frac12{\mathbb{Z}}_{\geq 0}$.}
\]

In what follows we will assume that  ${\mathfrak{C}}_{\lambda}$ is considered for
$\lambda\in{\mathbb C}$ because all the proofs hold for such $\lambda$,
not only for $\lambda\in{\mathbb{Z}}_{\geq 0}$.

Formulas (\ref{5.1.3}) easily imply that by setting
\[
H^t=H, \qquad U^t=-V, \qquad V^t=U
\]
we determine an antiautomorphism of superalgebra ${\mathfrak{C}}_{\lambda}$, i.e.,
$\left(u^t\right)^t=(-1)^{p(u)}u$ for any $u\in {\mathfrak{C}}_{\lambda}$.

\medskip

\noindent
{\bf Lemma 5.1.} {\it On ${\mathfrak{C}}_{\lambda}$, there exists a unique, up to a
scalar multiple, invariant linear functional $L$, such that $L(H)=0$
and $L\left(u^t\right)=L(u)$ for any $u\in {\mathfrak{C}}_{\lambda}$.  The functional $L$ is
uniquely determined by its restriction onto ${\mathbb C}[H]$ and its
generating function is}
\[
\varphi_{L}(t)=c\frac{e^{(\lambda +1)t}+e^{-\lambda t}}{e^t+1}.
\]

\noindent
{\bf Proof.} We have
\[
{\mathfrak{C}}_{\lambda}= \left\{\begin{array}{ll}
\mathop{\oplus}\limits_{i\in{\mathbb{Z}}}({\mathfrak{C}}_{\lambda})_{i}, &\text{if$\ \
\lambda\not\in\frac12{\mathbb{Z}}_{\geq 0}$},
\vspace{1mm}\\
\mathop{\oplus}\limits_{|i|\leq 2\lambda}({\mathfrak{C}}_{\lambda})_{i}, \quad &
\text{otherwise}, \end{array}\right.
\]
where $({\mathfrak{C}}_{\lambda})_{i}=\{u\in {\mathfrak{C}}_{\lambda}\mid [H, u]=iu\}$ for
$i\in{\mathbb{Z}}$. Hence, $L(({\mathfrak{C}}_{\lambda})_{i})=0$ unless $i=0$. Further
on,
\[
{}[U, fG]=\frac12\left (\lambda(\lambda+1)-H^2-H\right )f(H+1)+\frac12
\left (\lambda(\lambda+1)-H^2+H\right )f(H)
\]
and $L([U, fG])=0$. Therefore, as in Lemma 3.4, we obtain
\[
\left (\lambda(\lambda+1)-\left(\frac{d}{dt}\right)^2+\frac{d}{dt}\right )
\left(e^t+1\right)\varphi_{L}(t)=0
\]
and
\[
\varphi_{L}(t)=\left\{\begin{array}{ll} \displaystyle \frac{c_{1}e^{(\lambda +1)t}+
c_{2}e^{-\lambda t}}{1+e^{t}}, \quad &\text{for $\lambda\neq -\frac12$},\vspace{2mm}\\
\displaystyle \frac{(c_{1}+c_{2}t)e^{\frac12t}}{1+e^{t}}, &\text{otherwise.}\end{array}\right.
\]
The condition $L(H)=0$ implies that $c_{1}=c_{2}$ for $\lambda\neq
-\frac12$ and $c_{2}=0$ otherwise. This proves the uniqueness.

Let us prove the existence. By Lemma 3.1 we have
\[
[{\mathfrak{C}}_{\lambda},
{\mathfrak{C}}_{\lambda}]=[U, {\mathfrak{C}}_{\lambda}]+[V,
{\mathfrak{C}}_{\lambda}]
\]
and
\[\arraycolsep=0em
\begin{array}{l}
    [{\mathfrak{C}}_{\lambda},
{\mathfrak{C}}_{\lambda}]\cap{\mathbb C}[H]=([U, {\mathfrak{C}}_{\lambda}]+[V,
{\mathfrak{C}}_{\lambda}]\cap{\mathbb C}[H]=[U, ({\mathfrak{C}}_{\lambda})_1]\cap{\mathbb C}[H]
\vspace{2mm}\\
\qquad {} =\mbox{Span}([U,
fV]\mid f\in{\mathbb C}[H])=\mbox{Span}(f(H+1)Uv+f(H)VU\mid f\in{\mathbb C}[H]).
\end{array}
\]
Hence, $\deg [U, fV]=\deg f+2$. Hence, $1, H\not \in [{\mathfrak{C}}_{\lambda},
{\mathfrak{C}}_{\lambda}]$ and ${\mathfrak{C}}_{\lambda}/[{\mathfrak{C}}_{\lambda},
{\mathfrak{C}}_{\lambda}]$ is the linear span of the images of~1 and $H$. This
proves the existence of~$L$.

Let $\tilde L(u)=L\left(u^t\right)=0$. Then $\tilde L$ is invariant and $\tilde
L(H)=L(H)=0$. The uniqueness implies that $\tilde L=L$.\hfill\rule{3mm}{3mm}

\medskip

\noindent
{\bf Lemma 5.2.} {\it Let $f\in{\mathbb C}[H]$. We have

\vspace{-2mm}

\begin{enumerate}
\topsep0mm
\partopsep0mm
\parsep0mm
\itemsep0mm

\item[i)] if $f(2k-H)=-f(H)$, then $L\left(f(H)V^{2k}U^{2k}\right)=0$;

\item[ii)] if $f(2k+1-H)=f(H)$, then $L\left(f(H)V^{2k+1}U^{2k+1}\right)=0$.
\end{enumerate}}

\noindent
{\bf Proof.} i) Determine an automorphism of ${\mathfrak{C}}_{\lambda}$ by
setting
\[
\theta(H)=-H,\qquad \theta(U)=V,\qquad \theta(V)=U
\]
and set
\[
g(H)=f(H)V^{2k}U^{2k}+U^{2k}f(H)V^{2k}.
\]
Then
\[\arraycolsep=0em
\begin{array}{l}
    g(-H)=\theta(g(H))=\theta\left(f(H)V^{2k}U^{2k}+U^{2k}f(H)V^{2k}\right)
\vspace{2mm}\\
\phantom{g(-H)}  {}=
f(-H)V^{2k}U^{2k}+U^{2k}f(-H)V^{2k}
\vspace{2mm}\\
\phantom{g(-H)} {}=U^{2k}f(2k-H)V^{2k}+
f(2k-H)V^{2k}U^{2k}=-g(H).
\end{array}
\]
In other words, $g(H)$ is an odd polynomial. By Lemma 5.1 the
generating function for $L$ is an even one, so $L(g(H))=0$.

Further on,
\[
L\left(f(H)V^{2k}U^{2k}\right)=L\left(\frac12\left[f(H)V^{2k},
U^{2k}\right]\!+\!\frac12\left(f(H)V^{2k}U^{2k}\!+\!U^{2k}f(H)V^{2k}\right)\!\right)=0.
\]

ii) Set
\[
g(H)=f(H)V^{2k+1}U^{2k+1}-U^{2k+1}f(H)V^{2k+1}.
\]
Then
\[\arraycolsep=0em
\begin{array}{l}
    g(-H)=\theta(g(H))=\theta\left(f(H)V^{2k+1}U^{2k+1}-U^{2k+1}f(H)V^{2k+1}\right)
\vspace{2mm}\\
\phantom{g(-H)} {}=f(-H)V^{2k+1}U^{2k+1}-U^{2k+1}f(-H)V^{2k+1}
\vspace{2mm}\\
\phantom{g(-H)} {}= U^{2k+1}f(2k+1-H)V^{2k+1}+
f(2k+1-H)V^{2k+1}U^{2k+1}
\vspace{2mm}\\
\phantom{g(-H)} {}= V^{2k+1}f(H)U^{2k+1}-f(H)V^{2k+1}U^{2k+1}=-g(H).
\end{array}
\]
Therefore, $L(g(H))=0$. Further on,
\[
\arraycolsep=0em
\begin{array}{l}
\displaystyle f(H)V^{2k+1}U^{2k+1}=\frac12\left[f(H)V^{2k+1}, U^{2k+1}\right]
\vspace{3mm}\\
\displaystyle \phantom{f(H)V^{2k+1}U^{2k+1}={}} {}+
\frac12\left(f(H)V^{2k+1}U^{2k+1}-U^{2k+1}f(H)V^{2k+1}\right);
\end{array}
\]
hence, $L\left(f(H)V^{2k+1}U^{2k+1}\right)=0$.
\hfill\rule{3mm}{3mm}

\medskip

\noindent
{\bf Lemma 5.3.} {\it Let $f, g\in{\mathbb C}[H]$, $\varepsilon=\pm 1$,
$i\in{\mathbb{Z}}_{\geq 0}$. We have:

\vspace{-2mm}

\begin{enumerate}
\topsep0mm
\partopsep0mm
\parsep0mm
\itemsep0mm

\item[i)] If $f(H)V^{2i+1}=\left[U, gV^{2i+2}\right]$, then the condition
$\varepsilon g(H)=g(2i+2-H)$ is equivalent to the
condition $f(2i+1-H)=-\varepsilon f(H)$ and $\deg f=\deg g+1$.

\item[ii)] If $f(H)V^{2i}=\left[U, gV^{2i+1}\right]$, then the condition $\varepsilon
g(H)=g(2i+1-H)$ is equivalent to the
condition $f(2i-H)=\varepsilon f(H)$ and $\deg f=\deg g+2$.

\item[iii)] If $f(H)V^{2i+1}=\left\{U, gV^{2i+2}\right\}$, then the condition
$\varepsilon g(H)=g(2i+2-H)$ is equivalent to the
condition $f(2i+1-H)=\varepsilon f(H)$ and $\deg f=\deg g+2$.

\item[iv)] If $f(H)V^{2i}=\left\{U, gV^{2i+1}\right\}$, then the condition
$\varepsilon g(H)=g(2i+1-H)$ is equivalent to the
condition $f(2i-H)=-\varepsilon f(H)$ and $\deg f=\deg g+1$.
\end{enumerate}}

\noindent
{\bf Proof.} i) The condition $f(H)V^{2i+1}=[U, gV^{2i+2}]$ is equivalent
to the equation
\begin{equation}
f(H)=T(H+1)g(H+1)-T(H-2i-1)g(H), \label{5.4.1}
\end{equation}
where $T(H)=\frac12(\lambda+H)(\lambda+1-H)$. This implies
\[
\arraycolsep=0em
\begin{array}{l}
    f(2i+1-H)=T(2i+2-H+1)g(2i+2-H)-T(-H)g(2i+1-H)
\vspace{2mm}\\
\phantom{ f(2i+1-H)} {}= \varepsilon T(H-2i-1)g(H)-\varepsilon T(H+1)g(H+1)=\varepsilon f(H).
\end{array}
\]
Conversely, let $f(2i+1-H)=-\varepsilon f(H)$. Set $\psi(H)=g(H)-\varepsilon
g(2i+2-H)$. Then (\ref{5.4.1}) implies
\begin{equation}
\psi(H+1)T(H+1)-\psi(H)T(H-2i-1)=0.\label{5.4.2}
\end{equation}
If $i\geq 0$, then polynomial $T(H+1)$ has a root $\alpha$ such that
$\alpha-k$ is not a root of $T(H-2i-1)$ for any $k\in{\mathbb{Z}}_{\geq 0}$.
Select this root $\alpha$.  Then equation (\ref{5.4.2}) implies
$\psi(\alpha)=0$, but then $\psi(\alpha-1)T(\alpha-1-2i-1)=0$ and,
therefore, $\psi(\alpha-1)=0$.  So $\psi(\alpha-k)=0$ for any
$k\in{\mathbb{Z}}_{\geq 0}$.  Thus, $\psi=0$.

Headings ii)--iv) are similarly proved.\hfill\rule{3mm}{3mm}

\medskip

{\bf Proof of heading 2 of Theorem 2.4.} Let $l=2k$.  Let us
prove by induction on $k$ that $f_{2k, 2i}$ is an orthogonal
polynomials of degree $2k$ with respect to the form $\langle \cdot ,
\cdot\rangle_{2i}$. Let
\begin{equation}
{\mathbb C}[H]_{+}^i=\{f\in{\mathbb C}[H]\mid 
f(2i-H)=f(H)\}\label{5.4.3_{+}€}
\end{equation}
and
\begin{equation}
{\mathbb C}[H]_{-}^i=\{f\in{\mathbb C}[H]\mid f(2i-H)=-f(H)\}.
\label{5.4.3_{-}€}
\end{equation}
Then Lemma 5.2 implies that the spaces ${\mathbb C}[H]_{+}^i$
and ${\mathbb C}[H]_{-}^i$ are orthogonal
with respect to the form $\langle \cdot ,
\cdot\rangle_{2i}$. So it suffices to prove that $\langle f_{2k, 2i},
g\rangle_{2i}=0$ for $g\in{\mathbb C}[H]_{+}^i$ and
$\deg g<2k$.

Let us induct on $k$. If $k=1$, then
\[\arraycolsep=0em
\begin{array}{l}
\langle f_{2k, 2i}, 1\rangle_{2i}=\langle \left[U, V^{2i+1}\right],
V^{2i}\rangle\stackrel{\text{Lemma 3.2}}{=}-\langle V^{2i+1},
\left[U, V^{2i}\right]\rangle
\vspace{2mm}\\
\phantom{\langle f_{2k, 2i}, 1\rangle_{2i}} {}=-2\langle V^{2i+1}, V^{2i+1}\rangle\stackrel{\text{Lemma
5.2}}{=}0.
\end{array}
\]
Let $k>1$; then equations (\ref{2.8.3}) imply that
\begin{equation}
f_{2k, 2i}V^{2i}=\left[U, \frac{1}{H-\left(i+\frac12\right)}\left[f_{2k-2, 2i+2},
V^{2i+2}\right]\right].\label{5.5.1}
\end{equation}
This equation and Lemma 5.3 imply that $\deg f_{2k, 2i}=2k$.
Let $g\in{\mathbb C}[H]_{+}^i$
and $\deg g\leq 2k-2$. Then
\[
\arraycolsep=0em
\begin{array}{l}
    \displaystyle \langle f_{2k, 2i}, g\rangle_{2i}=\langle f_{2k, 2i}V^{2i},
gV^{2i}\rangle=\langle \left[U, \frac{1}{H-\left(i+\frac12\right)}\left[U, f_{2k-2, 2i+2}
V^{2i+2}\right]\right], gV^{2i}\rangle
\vspace{3mm}\\
\displaystyle \phantom{\langle f_{2k, 2i}, g\rangle_{2i}} {}=
-\langle \frac{1}{H-\left(i+\frac12\right)}\left[U, f_{2k-2, 2i+2}
V^{2i+2}\right], \left[V, gV^{2i}\right]\rangle
\end{array}
\]
\[\arraycolsep=0em
\begin{array}{l}
\displaystyle \phantom{\langle f_{2k, 2i}, g\rangle_{2i}} {}
= \langle \frac{1}{H-\left(i+\frac12\right)}\left[U, f_{2k-2, 2i+2}
V^{2i+2}\right], \nabla gV^{2i}\rangle
\vspace{3mm}\\
\displaystyle \phantom{\langle f_{2k, 2i}, g\rangle_{2i}} {} =
\langle \left[U, f_{2k-2, 2i+2}
V^{2i+2}\right], \frac{\nabla g}{H-\left(i+\frac12\right)}V^{2i+1}\rangle
\vspace{3mm}\\
\displaystyle \phantom{\langle f_{2k, 2i}, g\rangle_{2i}} {} =
\langle f_{2k-2, 2i+2}
V^{2i+2}, \left[V, \frac{\nabla g}{H-\left(i+\frac12\right)}V^{2i+1}\right]\rangle
\vspace{3mm}\\
\displaystyle \phantom{\langle f_{2k, 2i}, g\rangle_{2i}} {}=
\langle f_{2k-2, 2i+2}
V^{2i+2}, \frac{2\nabla g}{H-\left(i+\frac12\right)}-\nabla\left(\frac{\nabla
g}{H-(i+\frac12)}\right) \rangle_{2i+2}=0
\end{array}
\]
since $\deg \left(\frac{2\nabla g}{H-\left(i+\frac12\right)}-\nabla\left(\frac{\nabla
g}{H-\left(i+\frac12\right)}\right)\right)\leq 2k-3$.

Proof of orthogonality of polynomials $f_{2k+1, 2i}$ is similarly
performed with appellation to Lemma 3.3.

{\bf Proof of headings 3 and 4 of Theorem 2.4.} Define an
operator $D$ by setting
\[
Df\cdot V^{2i}=\left[V, \frac{\nabla g}{H-\left(i+\frac12\right)}\left[V, fV^{2i}\right]\right],
\]
where $f$ satisfies $f(2i-H)=f(H)$.  Let us show that $D$ is well
defined.  Indeed,
\[
(\nabla f)\left(i+\frac12\right)=f\left(i+\frac12\right)-f\left(i-\frac12\right)
=f\left(2i-\left(i-\frac12\right)\right)-f\left(i-\frac12\right)=0.
\]
Set $\varphi(H)=\displaystyle\frac{\nabla f}{H-(i+\frac12)}$. We have
\[
\varphi(2i+1-H)=\frac{f(2i+1-H)-f(2i-H)}{(2i+1-H)-\left(i+\frac12\right)}=
\frac{f(H-1)-f(H)}{i+\frac12-H}=\varphi(H).
\]

Therefore, $(Df)(2i-H)=(Df)(H)$ by Lemma 5.3.

Let us show now that $D$ is selfadjoint with respect to the form
$\langle \cdot , \cdot\rangle_{2i}$.  Indeed, let ${\mathbb
C}[H]_{+}^i$ and ${\mathbb C}[H]_{-}^i$ be as in (5.6)--(5.7).  We
have shown that $D$ sends ${\mathbb C}[H]_{+}^i$ into itself.  Let $f,
g\in {\mathbb C}[H]_{+}^i$; then
\[
\arraycolsep=0em
\begin{array}{l}
 \displaystyle    \langle Df,  g\rangle_{2i}=\langle Df\cdot V^{2i},  gV^{2i}\rangle=
\langle  \frac{1}{H-\left(i+\frac12\right)}\left[V, fV^{2i}\right],  \left[V, gV^{2i}\right]\rangle
\vspace{3mm}\\
\displaystyle \phantom{\langle Df,  g\rangle_{2i}} {}=
-\langle  \left[V, fV^{2i}\right],  \frac{1}{H-\left(i+\frac12\right)}\left[V, gV^{2i}\right]\rangle
\vspace{3mm}\\
\displaystyle
\phantom{\langle Df,  g\rangle_{2i}} {} =
\langle  fV^{2i},  \left[U, \frac{1}{H-\left(i+\frac12\right)}\left[V, gV^{2i}\right]\right]\rangle=
\langle f,  Dg\rangle_{2i}.
\end{array}
\]
Hence, $\langle Df_{2k, 2i}, g\rangle_{2i}=\langle f_{2k, 2i},
Dg\rangle_{2i}$ if $g\in {\mathbb C}[H]_{+}^i$ and $\deg g< 2k$.  The uniqueness of the
orthogonal polynomial of given degree implies that $Df_{2k,
2i}=\alpha_{k}f_{2k, 2i}$.  Furthermore,
\[
\arraycolsep=0em
\begin{array}{l}
\displaystyle
    Df\cdot V^{2i}=\left[U, \frac{1}{H-\left(i+\frac12\right)}\left[V, fV^{2i}\right]\right]=
\left[U, -\frac{\nabla f}{H-\left(i+\frac12\right)}V^{2i+1}\right]
\vspace{3mm}\\
\displaystyle
\phantom{Df\cdot V^{2i}} {}=-U\frac{\nabla
f}{H-\left(i+\frac12\right)}V^{2i+1}-\frac{\nabla f}{H-\left(i+\frac12\right)}V^{2i+1}U
\vspace{3mm}\\
\displaystyle
\phantom{Df\cdot V^{2i}} {} =
-\frac{\nabla f}{H-\left(i+\frac12\right)}UV\cdot V^{2i}-\frac{\nabla
f}{H-\left(i+\frac12\right)}V^{2i}VU
\vspace{3mm}\\
\displaystyle
\phantom{Df\cdot V^{2i}} {} =
\left(\frac12\cdot \frac{(H-\lambda)(H+\lambda+1)}{2H-i+\frac12}\Delta f+
\frac{(H-2i-\lambda+1)(H-2i+\lambda+1)}{2H-i+\frac12}\right )\Delta f.
\end{array}
\]
In other words,
\[
Df=\left(\frac{(H-\lambda)(H+\lambda+1)}{2H-2i+1}+
\frac{(H-2i-\lambda+1)(H-2i+\lambda+1)}{2H-2i+1}\right )\Delta f.
\]
By calculating the leading coefficient of $Df$ leads us to the
equation of heading~3.

Equation of heading 4 is similarly obtained by considering operator
\[
Df\cdot V^{2i}=\left\{U, \frac{1}{H-\left(i+\frac12\right)}\left\{V, fV^{2i}\right\}\right\}.
\]

{\bf Proof of headings 5--7 of Theorem 2.4.} Statements of
heading~5 follow from the study of automorphism $\theta$ given by
formulas
\[
\theta(H)=-H, \qquad \theta(U)=V, \qquad \theta(V)=U
\]
in the same way as in heading 2 of Theorem 2.4.

To prove statements of heading 6, consider the following subspaces
of ${\mathbb C}[H]$:
\[
\arraycolsep=0em
\begin{array}{l}
    {\mathbb C}[H]_{+}^{i+\frac12}=\{f\in{\mathbb C}[H]\mid f(2i+1-H)=f(H)\},
\vspace{2mm}\\
{\mathbb C}[H]_{-}^{i+\frac12}=\{f\in{\mathbb C}[H]\mid f(2i+1-H)=-f(H)\}.
\end{array}
\]
By Lemma 5.2 these subspaces are completely isotropic
with respect to the form $\langle \cdot ,
\cdot\rangle_{2i+1}$.  Let us show that the form $\langle f_{2k, 2i+1} ,
f_{2k+1, 2i+1}\rangle_{2i}\neq 0$, while the remaining scalar
products vanish. Indeed, equations (\ref{2.8.3})--(\ref{2.8.6}) imply
\begin{equation}
f_{2k, 2i+1}V^{2i+1}=\frac{1}{H-\left(i+\frac12\right)}\left[V, f_{2k, 2i+2}V^{2i+2}\right]
\label{5.7.1}
\end{equation}
and
\begin{equation}
f_{2l+1, 2i+1}\left[V, f_{2l, 2i+2}V^{2i+2}\right].
\label{5.7.2}
\end{equation}
Hence,
\[\arraycolsep=0em
\begin{array}{l}
\displaystyle     \langle f_{2k, 2i+1}, f_{2l+1, 2i+1}\rangle_{2i+1}=\langle
\frac{1}{H-\left(i+\frac12\right)}\left[U, f_{2k, 2i+2}V^{2i+2}\right],\left[U, f_{2l,
2i+2}V^{2i+2}\right]\rangle
\vspace{3mm}\\
\displaystyle \qquad {}=
\langle \left[U, f_{2k, 2i+2}V^{2i+2}\right], \frac{1}{H-\left(i+\frac12\right)}\left[U, f_{2l,
2i+2}V^{2i+2}\right]\rangle
\vspace{3mm}\\
\displaystyle \qquad {} =
\left\langle f_{2k, 2i+2}V^{2i+2}, \left[V, \frac{1}{H-\left(i+\frac12\right)}\left[U, f_{2l,
2i+2}V^{2i+2}] \right.\right]\right\rangle
\vspace{3mm}\\
\displaystyle \qquad {} =
\langle f_{2k, 2i+2}V^{2i+2}, \tilde D f_{2l,
2i+2}V^{2i+2}\rangle=\langle f_{2k, 2i+2}, \tilde D f_{2l,
2i+2}\rangle_{2i+2},
\end{array}
\]
where
\[
\tilde D f\cdot V^{2i+2}=\left[V, \frac{1}{H-\left(i+\frac12\right)}\left[U, f_{2l,
2i+2}V^{2i+2}\right]\right]
\]
and is well defined thanks to Lemma 5.3i).

It is easy to show that $\tilde D$ is selfadjoint and $f_{2k, 2i+2}$
is its eigenfunction corresponding to a nonzero eigenvalue
$\alpha_{l}$. So
\[
\langle f_{2k, 2i+1}, f_{2l+1,
2i+1}\rangle_{2i+1}=\alpha_{l}\langle f_{2k, 2i+2}, f_{2l,
2i+2}\rangle_{2i+2}.
\]
This proves statement of heading 6.

Statements of heading 7 are proved similar to those of heading~5.

{\bf Proof of headings 1, 8, 9 of Theorem 2.4.} Proof of heading 1 is
similar to that of heading 1 of Theorem 2.2.  To prove heading 8,
consider the operator
\[
\tilde D f\cdot V^{2i+1}=\left[U, \left[V, fV^{2i+1}\right]\right]\quad \text{for}\ \ f\in
{\mathbb C}[H]_{+}^{i+\frac12}.
\]
It is easy to verify that $\tilde D f\in {\mathbb C}[H]_{-}^{i+\frac12}$ and,
the other way round, if $f\in {\mathbb C}[H]_{-}^{i+\frac12}$, then $\tilde D
f\in {\mathbb C}[H]_{+}^{i+\frac12}$.  Moreover, $\langle Df, g\rangle_{2i+1}=
\langle f, Dg\rangle_{2i+1}$, i.e., $D$ is
selfadjoint and $\deg Df=\deg f+1$ if $\deg f>0$.  Hence,
\[
\langle Df_{2k, 2i+1}, f_{2l, 2i+1}\rangle=
\langle f_{2k, 2i+1}, Df_{2l, 2i+1}\rangle =0\quad \text{if}\ \ l<k.
\]
So $Df_{2k, 2i+1}=\alpha\cdot Df_{2k+1, 2i+1}$.  Having calculated
$\alpha$ we obtain the statement of heading~8.

Proof of heading 9 is similar to arguments from the proof of of 
headings 8 and 9 of Theorem 2.2.

\subsection*{Acknowledgements}

I am thankful to D~Leites for encouragement and help and to
ESI, Vienna, for hospitality and support.

\label{sergeev-lastpage}
\newpage


\begin{thebibliography}{99}
\small
\topsep0mm
\partopsep0mm
\parsep0mm
\itemsep0mm

\bibitem{sergeev:Di}
Dixmier J, Alg\`ebres Envellopentes, Gautier-Villars, Paris,
1974; Enveloping Algebras, AMS, 1996.

\bibitem{sergeev:F}
Feigin B~L, The Lie Algebras ${\mathfrak{gl}} (\lambda)$ and Cohomologies of
Lie Algebra of Differential Operators, {\it Russian Math.  Surveys}
{\bf 43}, Nr.2 (1988),  157--158.

\bibitem{sergeev:GL}
Grozman P and Leites D, Lie Superalgebras of Supermatrices of Complex
Size and Integrable Dynamical Systems,
in Complex Analysis and Related Topics, Proc.  of the
International Symposium, Editors N~Vasilevsky et.~al.,
Cuernavaca, Mexica, November 18--22, 1996,
Birkhauser, 1999, 73--105.

\bibitem{sergeev:KLP}
Kwon K~H, Lee D~W and  Park S~B, New Characterization of Discrete
Classical Orthogonal Polynomials, {\it J.~Approx. Theory} {\bf 89}
(1977), 156--171.

\bibitem{sergeev:KS}
Koekoek R and Swarttouw R, The Askey-Scheme of Hypergeometric
Orthogonal Polynomials and its $q$-Analogue, math.CA/9602214.

\bibitem{sergeev:LS}
 Leites D and Sergeev A, Orthogonal Polynomials of Discrete Variable
 and Lie Algebras of Complex Size Matrices, {\it Teoret.  Mat.  Fiz.},
 (Russian) V. 123 (2000), no.  2, 205--236 translation in {\it
 Theoret.  and Math.  Physics} V. 123 (2000), no.  2, 582--608;
 preprinted in Procedings of M~Saveliev Memorial Conference, MPI,
 Bonn, February, 1999, MPI-1999-36, 49--70.

\bibitem{sergeev:M}
Marcellan F, Complex Path Integral Representation for Semiclassical
Linear Functionals, {\it J.~Approx. Theory} {\bf 94} (1998), 107--127.

\bibitem{sergeev:MBP}
Marcellan F, Branquinho A and  Petronillo T, Classical Orthogonal
Polynomials: a Functional Approach, {\it Acta Applicandae Math.}
{\bf 34} (1994), 
283--303.

\bibitem{sergeev:MK}
Mingarelli A and Krall A, Jacobi-Type Polynomials under an Indefinite
Inner Product, {\it Proc. Roy. Soc. Edinburgh, Sect. A} {\bf 90}
(1981), 147--153.

\bibitem{sergeev:NSU}
Nikiforov A~F, Suslov S~K and Uvarov V~B, Classical Orthogonal
Polynomials of a Discrete Variable, Springer Series in Computational
Physics, Springer-Verlag - Berlin, 1991 (Translated from the Russian).

\bibitem{sergeev:Mo}
Montgomery S, Constructing Simple Lie Superalgebras from Associative
Graded Algebras, {\it J.~Algebra} {\bf 1 95}, Nr.2 (1997), 558--579.

\bibitem{sergeev:Pi}
Pinczon G, The Enveloping Algebra of the Lie Superalgebra
${\mathfrak{osp}}(1|2)$,
{\it J. Algebra} {\bf 132} (1990), 219--242.

\end{thebibliography}
\end{document}